\documentclass{article}
\usepackage{amsfonts,amsmath,amssymb}
\newtheorem{theorem}{Theorem}[section]
\newtheorem{lemma}[theorem]{Lemma}

\numberwithin{equation}{section}
\begin{document}
\voffset=-1 true cm \setcounter{page} {1}
\title{ Blow-up solutions for a Kirchhoff type elliptic equation with trapping potential  }
\author{Helin Guo$^a$, Yimin Zhang$^b$ and Huan-Song Zhou$^b$ \thanks{Corresponding
author.\newline Email address: H.L.Guo: qfguohelin@126.com;  Y.M.Zhang: zhangyimin@whut.edu.cn; H.S.Zhou: hszhou@whut.edu.cn. \hfill
\newline This work was supported by  NFSC Grants 11471331, 11501555 and 11471330).}\\
\small $^a$University of Chinese Academy of Sciences and \\
\small Wuhan Institute of Physics and Mathematics, CAS, Wuhan 430071, China\\
\small $^b$Center for Mathematical Sciences and Department of Mathematics, \\
\small Wuhan University of Technology, Wuhan 430070, China
}
\date{}
\maketitle
{\bf Abstract}\quad We study a Kirchhoff type elliptic equation with trapping potential. The existence and blow-up behavior
of solutions with normalized $L^{2}$-norm for this equation are discussed. 
\section{Introduction}
\quad\quad In this paper, we are interested in the existence and blow-up behavior of solutions with normalized $L^{2}$-norm (i.e., $L^{2}$-norm equals 1) for the
following type elliptic problem
\begin{equation}\label{eq1.1}
-\left(a+b\int_{\mathbb R^{N}}|\nabla u|^{2}dx\right)\triangle u+V(x)u=\beta|u|^{p}u+\lambda u,\quad \ x\in \mathbb{R}^{N},
\end{equation}
where $1\leq N\leq4$, $\lambda\in \mathbb{R}\setminus\{0\}$, $a$ and $b$ are positive constants, $\beta>0$, $V(x)\in C(\mathbb R^{N}, \mathbb R^+),\ p\in(0,2^{\ast}-2)$ with $2^{\ast}=\frac{2N}{N-2}$ if $N=3,4$ and $2^{\ast}=+\infty$ if $N=1,2$. \eqref{eq1.1} is a steady-state equation of certain type of generalized Kirchhoff equation, which is usually called a Kirchhoff type elliptic equation. Classical Kirchhoff equation was proposed in \cite{K}, which is essentially a modified one dimensional wave equation and can be used to give a more accurate description on the transversal oscillations of a stretched string, see e.g.,  \cite{ap,as,K} for more backgrounds and related results on the classical Kirchhoff equations.

Over the last decade, much attention has been paid to the Kirchhoff type elliptic equation \eqref{eq1.1}, for examples, when $V(x)$ is a nonnegative constant, the existence of radial solutions of \eqref{eq1.1} was proved in \cite{JW} for $p>2$. He and Zou in \cite{HZ} studied \eqref{eq1.1} with $\lambda=0$ and $p\in (2,4)$, in which the mountain pass theorem and the Nehari manifold were directly used to obtain a positive ground state solution to \eqref{eq1.1}. For the case $p \in (0,1]$, He and Li \cite{HL} obtained a positive ground state solution to \eqref{eq1.1} by constructing a special Palais-Smale sequence. Problem \eqref{eq1.1} with $\lambda=0$ and $p\in (1,4)$ was discussed in \cite{LY} and  a positive ground state solution was found by solving a constrained minimization over a Nehari-Pohozaev type manifold. Note that when $b=0$, problem \eqref{eq1.1} is related to the so-called Gross-Pitaevskii equation which
arises in the study of Bose-Einstein condensates, see e.g.,\cite{DG,F}. In the case of $b=0$, the solutions with normalized $L^{2}$-norm, i.e., $L^{2}$-norm is equal to
1, have special interest in physics \cite{GS,Z} and the existence of this kind of solutions
has been established in \cite{GS,Z} for \eqref{eq1.1} with $N=2$ and $p=2$. Particularly, some detailed analysis on the asymptotic behavior of this kind of solutions were also discussed in \cite{GS} as $\beta$ approaches a critical value. Motivated by \cite{GS}, the authors of \cite{GZZ} studied the behavior of normalized solutions of \eqref{eq1.1} for fixed $\beta$ and $b=0$, but $p\nearrow2$. Some more generalizations on the results of \cite{GZZ} can be founded in \cite{Z1, ZZ}.
Recently, also inspired by \cite{GS}, in papers \cite{Y1,Y2,ZZ1} the authors studied the existence of solution $u$ with $\|u\|_{L^{2}}=c$ ($c>0$ is a constant)
for the Kirchhoff type elliptic problem \eqref{eq1.1} with $\beta=1$ and $V(x)\equiv0$, and their results show that the existence of this kind solutions depend heavily on the constant $c$.
The main aim of this paper is to extend the results of \cite{GS,GZZ,Z1} on problem \eqref{eq1.1} with $b=0$ (local case) to the nonlocal case (i.e. $b\neq0$), that is, we are interested in the existence and asymptotic behavior of solutions with normalized $L^{2}$-norm for problem \eqref{eq1.1} when $b\neq0$ and $p$ approaches
the critical value $p^{\ast}\triangleq\frac{8}{N}$. For this purpose, we consider the following constrained minimization problem:
\begin{equation}\label{eq1.3}
d_{\beta}(p)=\inf_{u\in S_{1}}E_{p}^{\beta}(u)
\end{equation}
where
\begin{equation}\label{eq1.4}
E_{p}^{\beta}(u)=\frac{a}{2}\int_{\mathbb R^{N}}|\nabla u|^{2}dx+\frac{b}{4}\left(\int_{\mathbb R^{N}}|\nabla u|^{2}dx\right)^{2}+\frac{1}{2}\int_{\mathbb R^{N}}V(x)u^{2}dx-\frac{\beta}{p+2}\int_{\mathbb R^{N}}|u|^{p+2}dx,
\end{equation}
\begin{equation}\label{a1}
S_{1}=\left\{u\in \mathcal{H}:\int_{\mathbb R^{N}}|u|^{2}dx=1\right\}\quad\mbox{with}\quad\mathcal{H}\triangleq\Big\{u\in H^{1}(\mathbb R^{N}):\int_{\mathbb R^{N}}V(x)u^{2}dx<\infty\Big\}.
\end{equation}
Clearly, $\mathcal{H}=H^{1}(\mathbb R^{N})$ if $V(x)\equiv0$. Moreover, if $N\geq3$, then $p\in(0,\frac{4}{N-2}]$ is necessary to ensure that the functional given by \eqref{eq1.4} is well defined in $\mathcal{H}$.
On the other hand, for any fixed $u_{0}\in S_{1}$, it is easy to see that $u^\sigma_{0}(x)=\sigma^{\frac{N}{2}}u_{0}(\sigma x)\in S_{1}$ for any $\sigma>0,$ but
\begin{equation*}
E_{p}^{\beta}(u_{0}^\sigma)\rightarrow -\infty\quad\mbox{as}\ \sigma\rightarrow +\infty,\quad\mbox{if}\  p>\frac{8}{N},
\end{equation*}
this means that $d_{\beta}(p)=-\infty$ if $p>\frac{8}{N}$. Therefore, throughout the paper, we always assume that
\begin{equation*}
0<p\leq\frac{8}{N}, \ \mbox{and}\quad p\leq\frac{4}{N-2}\quad\mbox{if}\  N\geq3\Longrightarrow\ N\leq4.
\end{equation*}
These implies that $N$ can be $1,2,3$ or $4$. Note that $N<4$ is essentially required in \cite{Y1}, (see the derivation of $(2.6)$ in \cite{Y1}). When $V(x)\equiv0$ and $N < 4$, by almost the same tricks as that of \cite{Y1}, we know that \eqref{eq1.3} has no minimizers for all $\beta>0$ if $p\geq\frac{8}{N}$, but there exists $\beta^{\ast}>0$ such that \eqref{eq1.3} has a minimizer if and only if $\beta\in(\beta^{\ast},+\infty)$ as $p\in(0,\frac{4}{N}]$, or $\beta\in[\beta^{\ast},+\infty)$ as $p\in(\frac{4}{N},\frac{8}{N})$. However, using the methods of \cite{Y1}, $\beta^{\ast}$ can be calculated only for $p\in (0,\frac{4}{N}]$. In this paper, we successfully obtain the explicit expression of $\beta^*$ for all $p \in (0, \frac{8}{N})$ and  $1\leq N \leq4$. We mention that, when $V(x)\equiv0$ in \eqref{eq1.3}, although the existence of minimizers for \eqref{eq1.3} is essentially proved in \cite{Y1,Y2} ($N\leq 3$), here we provide a very simple and direct proof for the existence of minimizers of \eqref{eq1.3} with $V(x) \equiv 0$, and $N=4$ is also included, see section 2. Comparing to the case of $V(x) \equiv0$, the other aim of the paper is to know whether there is any new phenomena for problem \eqref{eq1.3} when $V(x)\not\equiv0$.
In fact, our results of this paper show that the situation of \eqref{eq1.3} with $V(x)\not\equiv0$ is totally different from that of $V(x)\equiv0$. Roughly speaking, we prove that \eqref{eq1.3} has always a minimizer for all $\beta>0$ when $p\in(0,\frac{8}{N})$, and there exists $\beta^{\ast}>0$ such that \eqref{eq1.3} with $p=\frac{8}{N}$ has a minimizer if and only if $\beta\in(0,\beta^{\ast}]$, and $\beta^{\ast}$ can be given explicitly, see our Theorems \ref{th1.1} and \ref{th1}. Moreover, we establish a detailed analysis on the asymptotic behavior of the minimizer of \eqref{eq1.3} as $p\nearrow p^{\ast}=\frac{8}{N}$, see our Theorem \ref{th2}.


For stating our results, we introduce the following semilinear elliptic equation:
\begin{equation}\label{eq1.9}
-\frac{Np}{4}\triangle u+\left(1+\frac{p}{4}(2-N)\right)u-u^{p+1}=0,\quad x\in\mathbb R^{N},\ 1\leq N\leq4,\ 0<p<2^{\ast}-2,
\end{equation}
it is well-known that this equation has a unique ({\it up to translations}) positive solution $\phi_{p}\in H^1(\mathbb{R}^N)$, which is radially symmetric and decays exponentially at infinity, see e.g., \textup{\cite{GNN,KJ,WM}}.

\begin{theorem}\label{th5}
If $V(x)\equiv0$, $p\in(0,p^*]$ with $p^*=\frac{8}{N}$ and $1\leq N\leq4$, let
\begin{equation}\label{a2}
\widetilde{\beta}_{p}=
\begin{cases}
0, \quad\quad\quad\quad\quad\quad\quad\quad\quad\quad\quad\quad\quad\quad\quad\  \mbox{if}\quad 0<p<\frac{4}{N},
\\\\a\|\phi_{p}\|_{L^{2}}^{\frac{4}{N}},\quad\quad\quad\quad\ \quad\quad\quad\quad\quad\quad\quad\quad\mbox{if}\quad p=\frac{4}{N},
\\\\2\|\phi_{p}\|_{L^{2}}^{p}\Big(\frac{2a}{8-Np}\Big)^{\frac{8-Np}{4}}\Big(\frac{b}{Np-4}\Big)^{\frac{Np-4}{4}},\quad\ \mbox{if}\quad\frac{4}{N}<p<p^*=\frac{8}{N},
\end{cases}
\end{equation}
and
\begin{equation}\label{w11}
\widetilde{\beta}_{p^{\ast}}=
\begin{cases}
\frac{b}{2}\|\phi_{p^{\ast}}\|_{L^{2}}^{\frac{8}{N}},\quad\mbox{if}\ N\leq3,\ p=p^{\ast},
\\\\bS^{2},\quad\quad\quad\ \mbox{if}\ N=4,\ p=p^{\ast},
\end{cases}
\end{equation}
where $\phi_{p}$ is the unique  positive solution of (\ref{eq1.9})
and
\begin{equation}\label{w12}
S=\inf_{u\in D^{1,2}(\mathbb R^{4})}\frac{\|\nabla u\|_{L^{2}}^{2}}{\|u\|_{L^{4}}^{2}}>0.
\end{equation}
Then,\\
$(\textup{i})$ For $p\in(0,\frac{8}{N})$, problem \eqref{eq1.3} has a minimizer if and only if
\begin{equation} \label{beta}
\beta>\widetilde{\beta}_{p}\ \mbox{with}\ 0<p\leq\frac{4}{N},\ \mbox{or}\
\beta\geq \widetilde{\beta}_{p}\ \mbox{with}\ \frac{4}{N}<p<\frac{8}{N}.
\end{equation}
Therefore,  under the condtions of \eqref{beta}, problem \eqref{eq1.1} has always a positive solution $u\in H^1(\mathbb{R}^N)$ with $\|u\|_{L^{2}}=1$ for some $\lambda<0$. \\
$(\textup{ii})$ For $p=p^*(=\frac{8}{N})$, \eqref{eq1.3} has no minimizers for any $\beta>0$.
\end{theorem}

Let
\begin{equation}\label{eq1.10}
\beta_{p}=\frac{b}{2}\|\phi_{p}\|_{L^{2}}^{p},
\end{equation}
where $\phi_{p}$ is the unique positive solution of \eqref{eq1.9}. Then, we have the following theorem.
\begin{theorem}\label{th1.1}
Let $0<p<p^{\ast}$ and $1\leq N\leq4$. If $V(x)$ satisfies
\begin{equation}\label{a3}
0\not\equiv V(x)\in C(\mathbb R^{N},\mathbb R^{+}),\quad \lim_{|x|\rightarrow\infty}V(x)=\infty\quad\mbox{and}\quad \inf_{x\in\mathbb R^{N}}V(x)=0.
\end{equation}
Then, for any fixed $\beta>0$, problem \eqref{eq1.3} has at least one minimizer.
\end{theorem}

 When $N=4$, since $p^{\ast}+2=2^{\ast}=4$ is the Sobolev critical exponent,  then we have to confront simultaneously the noncompactness problems caused by the unboundedness of the domain $\mathbb{R}^N$ and
 the Sobolev critical growth. In this case,  it is well-known that  even  the embedding of the radially symmetric space of $\mathcal{H}$ into $L^{2^*}(\mathbb R^{4})$ is not compact. 
 For these reasons, we can prove  the following results only for $N\leq3$.
\begin{theorem}\label{th1}
If $1\leq N\leq 3$ and $V(x)$ satisfies the condition \eqref{a3}
Then, for $p=p^{\ast}$, we have\\
$(\textup{i})$ $d_{\beta}(p^{\ast})>0$ and \eqref{eq1.3} has at least one minimizer if $0<\beta\leq \beta_{p^{\ast}}\triangleq\frac{b}{2}\|\phi_{p^{\ast}}\|_{L^{2}}^{p\ast}$.\\
$(\textup{ii})$ $d_{\beta}(p^{\ast})=-\infty$ and \eqref{eq1.3} has no minimizers if $\beta>\beta_{p^{\ast}}.$
\end{theorem}

Based on the above existence results, our following theorem gives some asymptotic properties of the minimizers of \eqref{eq1.3} as $p\nearrow p^{\ast}$.
\begin{theorem}\label{th3}
For any given $\beta \in (0,\beta_{p^{\ast}})$ and $1\leq N\leq3$, let $V(x)$ satisfy \eqref{a3} and let $u_{p}\in S_{1}$ be a nonnegative minimizer of problem \eqref{eq1.3} for each $p\in(0,p^{\ast})$. Then, there is a subsequence of $\{u_{p}\}$, still denoted by $\{u_{p}\}$, such that, for some $u_{0}\in\mathcal{H}$,
\begin{equation*}
d_{\beta}(p)\rightarrow d_{\beta}(p^{\ast})\quad\mbox{and}\quad u_{p}\rightarrow u_{0}\ \mbox{strongly}\ in\ \mathcal{H}\quad \mbox{as}\ p\nearrow p^{\ast}.
\end{equation*}
Moreover, $u_{0}\in S_{1}$ is a minimizer of $d_{\beta}(p^{\ast})$.
\end{theorem}

For any fixed $\beta>\beta_{p^{\ast}}$,  we know that there exists a positive constant  $m$ ( independent of $p$) such that
$\frac{\beta}{\beta_{p^{\ast}}}>m>1$. Then,  it follows from
 $\beta_{p}\rightarrow\beta_{p^{\ast}}$
as $p\nearrow p^{\ast}$  that
\begin{equation}\label{eq1.13}
\left(\frac{\beta p}{\beta_{p}p^{\ast}}\right)^{\frac{p^{\ast}}{p^{\ast}-p}}
\rightarrow+\infty\quad\mbox{as}\ p\nearrow p^{\ast}.
\end{equation}
Using the fact of \eqref{eq1.13}, we have the following theorem, which describes the concentration behavior of minimizers for \eqref{eq1.3} as $p\nearrow p^{\ast}$.
\begin{theorem}\label{th2}
For any fixed $\beta>\beta_{p^{\ast}}$ and $1\leq N\leq3$, let $V(x)$ satisfy \eqref{a3} and let $u_{p}$ be a nonnegative minimizer of \eqref{eq1.3} for each $p\in(0,p^{\ast})$. Then,
for any sequence of $\{u_{p}\}$ with $p\nearrow p^{\ast}$,  there exist $\{y_{\epsilon_{p}}\}\subset\mathbb R^{N}$ and
$y_{0}\in \mathbb R^{N}$ such that
\begin{equation*}
\lim_{p\nearrow p^{\ast}}\epsilon_{p}^{\frac{N}{2}}u_{p}(\epsilon_{p}x+\epsilon_{p}y_{\epsilon_{p}})=\frac{1}{\|\phi_{p^{\ast}}\|_{L^{2}}}\phi_{p^{\ast}}\left(|x-y_{0}|\right) \quad \mbox{in} \ H^{1}(\mathbb R^{N}),
\end{equation*}
where\\
\begin{equation*}
\epsilon_{p}=\left(\frac{\beta p}{\beta_{p}p^{\ast}}\right)^{-\frac{p^{\ast}}{4(p^{\ast}-p)}}\rightarrow 0.
\end{equation*}
Moreover, $\{y_{\epsilon_{p}}\}$ satisfies $\epsilon_{p}y_{\epsilon_{p}}\rightarrow z_{0}$ as $p\nearrow p^{\ast}$, and $z_{0}$ is a global minimal point
of $V(x)$, i.e., $V(z_{0})=0.$
\end{theorem}
\section{Existence and nonexistence for \eqref{eq1.3} with $V(x)\equiv0$}
In this section, we prove first Theorem \ref{th5}, and then establish some energy estimates which are required in next section. Before going to the proofs, we recall the following Gagliardo-Nirenberg inequality \cite{W2}
\begin{equation}\label{eq1.8}
\|u\|_{L^{p+2}}^{p+2}\leq \frac{p+2}{2\|\phi_{p}\|_{L^{2}}^{p}}\|\nabla u\|_{L^{2}}^{\frac{Np}{2}}\|u\|_{L^{2}}^{2+\frac{p}{2}(2-N)},\ N\geq1,\ 0<p<2^{\ast}-2,
\end{equation}
where $\phi_{p}$ is given in \eqref{a2}. Moreover, by \eqref{eq1.9} and the Pohozaev identity, we know that
\begin{equation}\label{eq1.11}
\int_{\mathbb R^{N}}\phi_{p}^{2}dx=\int_{\mathbb R^{N}}|\nabla\phi_{p}|^{2}dx,\quad
\int_{\mathbb R^{N}}\phi_{p}^{2}dx=\frac{2}{p+2}\int_{\mathbb R^{N}}|\phi_{p}|^{p+2}dx.
\end{equation}
When $V(x)\equiv0$, for the sake of simplicity, we rewrite \eqref{eq1.3} as follows:
\begin{equation}\label{eq1.7}
\widetilde{d}_{\beta}(p)=\inf_{u\in \widetilde{S}_{1}}\widetilde{E}_{p}^{\beta}(u),
\end{equation}
where $\widetilde{E}_{p}^{\beta}$ is given by
\begin{equation} \label{Ebeta}
\widetilde{E}_{p}^{\beta}(u)=\frac{a}{2}\int_{\mathbb R^{N}}|\nabla u|^{2}dx+\frac{b}{4}\left(\int_{\mathbb R^{N}}|\nabla u|^{2}dx\right)^{2}-\frac{\beta}{p+2}\int_{\mathbb R^{N}}|u|^{p+2}dx
\end{equation}
and
\begin{equation*}
\widetilde{S}_{1}=\left\{u\in H^{1}(\mathbb R^{N}):\int_{\mathbb R^{N}}|u|^{2}dx=1\right\}.
\end{equation*}
Clearly, $\widetilde{d}_{\beta}(p)$ is well defined and
\begin{equation}\label{eqabc}
\widetilde{d}_{\beta}(p) \leq 0 \ \mbox{for all}\ \beta>0, \text{ and } \widetilde{d}_{\beta}(p) = 0 \ \mbox{for all}\ \beta\leq 0.
\end{equation}
 In fact,
for any $\beta>0$ and $u\in \widetilde{S}_{1}$, using \eqref{eq1.8} and the definition of (\ref{Ebeta}) we see that
\begin{equation}\label{a4}
\widetilde{E}_{p}^{\beta}(u)\geq \frac{a}{2}\int_{\mathbb R^{N}}|\nabla u|^{2}dx + \frac{b}{4}\left(\int_{\mathbb R^{N}}|\nabla u|^{2}dx\right)^{2}
-\frac{\beta}{2\|\phi_{p}\|_{L^{2}}^{p}}\left(\int_{\mathbb R^{N}}|\nabla u|^{2}dx\right)^{\frac{Np}{4}},
\end{equation}
this implies that $\widetilde{E}_{p}^{\beta}$ is bounded from below on $\widetilde{S}_{1}$ since $0<p<\frac{8}{N}$ and $0<\frac{Np}{4}<2$, and (\ref{eq1.7}) is well defined. By \eqref{a4} , it is obvious that $\widetilde{d}_{\beta}(p) \geq 0 \ \mbox{for all}\ \beta\leq 0$. Moreover, taking $u\in \widetilde{S}_{1}$ and
letting $u_{t}(x)=t^{\frac{N}{2}}u(tx)(t>0)$, then $u_{t}\in \widetilde{S}_{1}$ and
\begin{equation}\label{eq2.2}
\widetilde{d}_{\beta}(p)\leq \widetilde{E}_{p}^{\beta}(u_{t})=\frac{at^{2}}{2}\int_{\mathbb R^{N}}|\nabla u|^{2}dx
+\frac{bt^{4}}{4}\left(\int_{\mathbb R^{N}}|\nabla u|^{2}dx\right)^{2}
-\frac{\beta t^{\frac{Np}{2}}}{p+2}\int_{\mathbb R^{N}}|u|^{p+2}dx\rightarrow0\quad \mbox{as}\ t\rightarrow0.
\end{equation}
Hence,  (\ref{eqabc}) is proved.

\begin{lemma}\label{lea}
Let $1\leq N\leq4$ and $p\in(0,\frac{8}{N})$. Then, $\widetilde{d}_{\beta}(p)<0$ if and only if $\beta>\widetilde{\beta}_{p}$, where $\widetilde{\beta}_{p}$ is defined by \eqref{a2}. Moreover,  $\widetilde{d}_{\beta}(p)=0$ for all $\beta \leq \widetilde{\beta}_{p}$ and $p\in (0,\frac{8}{N})$.
\end{lemma}
Proof. We prove this lemma by three cases: Case 1. $p\in(0,\frac{4}{N})$; Case 2. $p=\frac{4}{N}$ and Case 3. $p\in(\frac{4}{N},\frac{8}{N})$.\\
The first two cases can be proved by similar way to that of Lemma 2.3 in \cite{Y1}. But, for the third case, we have to use a new approach which allows us to get an explicit expression of $\widetilde{\beta}_{p}$ and to include $N=4$, these are impossible by following \cite{Y1}.


{\it Case 1:} $p\in(0,\frac{4}{N})$. In this case, $0<\frac{Np}{2}<2$ and $\widetilde{\beta}_{p}=0$ by \eqref{a2}, then  (\ref{eq2.2}) shows that $\widetilde{d}_{\beta}(p)<0$ for all $\beta>0=\widetilde{\beta}_{p}$, and $\widetilde{d}_{\beta}(p) = 0$ for all $ \beta\leq 0=\widetilde{\beta}_{p}$ by \eqref{eqabc}.

{\it Case 2:} $p=\frac{4}{N}$. In this case, $\widetilde{\beta}_{p}=\widetilde{\beta}_{4/N}$ by \eqref{a2}, and we have two different  situations.\\
$\bullet$ If $\beta  \leq \widetilde{\beta}_{p}$,  then $\widetilde{d}_{\beta}(p)=0$ for each $\beta \leq \widetilde{\beta}_{p}$. Indeed,
let $u\in \widetilde{S}_{1}$, it follows from\eqref{eq1.8} and \eqref{a4} that
\begin{equation*}
\widetilde{E}_{p}^{\beta}(u)\geq\frac{a(\widetilde{\beta}_{p}-\beta)}{2\widetilde{\beta}_{p}}\int_{\mathbb R^{N}}|\nabla u|^{2}dx
+\frac{b}{4}\left(\int_{\mathbb R^{N}}|\nabla u|^{2}dx\right)^{2}>0, \text{ note now that } Np/4=1 \text{ in } \eqref{a4},
\end{equation*}
hence, $\widetilde{d}_{\beta}(p)\geq0$, and then $\widetilde{d}_{\beta}(p) = 0$ since \eqref{eqabc}.\\
$\bullet$ If $\beta>\widetilde{\beta}_{p}$, then $\widetilde{d}_{\beta}(p)< 0$. In fact, let $\phi_{p}$ be given by \eqref{eq1.8} and  set
\begin{equation} \label{ut}
u_{t}=\frac{t^{\frac{N}{2}}\phi_{p}(tx)}{\|\phi_{p}\|_{L^{2}}}, \text{ for } t>0,
\end{equation}
then $u_{t}\in \widetilde{S}_{1}$
and it follows from \eqref{eq1.11} that
\begin{equation} \label{eq2.5}
\int_{\mathbb R^{N}}|\nabla u_{t}|^{2}dx
=\frac{t^{2}\int_{\mathbb R^{N}}|\nabla\phi_{p}|^{2}dx}{\|\phi_{p}\|_{L^{2}}^{2}}=t^{2},
\end{equation}
\begin{equation} \label{eq2.6}
\int_{\mathbb R^{N}}u_{t}^{p+2}dx
=\frac{t^\frac{Np}{2}\int_{\mathbb R^{N}}|\phi_{p}|^{p+2}dx}{\|\phi_{p}\|_{L^{2}}^{p+2}}
=\frac{(p+2)t^{\frac{Np}{2}}}{2\|\phi_{p}\|_{L^{2}}^{p}}.
\end{equation}
Note that, in this case,  $p=\frac{4}{N}$ and $\frac{Np}{2} =2$. Hence,
\begin{equation*}
\widetilde{E}_{p}^{\beta}(u_{t})=\left(\frac{a}{2}-\frac{\beta}{2\|\phi_{p}\|_{L^{2}}^{\frac{4}{N}}}\right)t^{2}+\frac{b}{4}t^{4}.
\end{equation*}
Therefore, $\widetilde{d}_{\beta}(p)\leq\inf\limits_{t>0}\widetilde{E}_{p}^{\beta}=-\frac{a^{2}}{4b}\left(\frac{\beta}{\widetilde{\beta}_{p}}-1\right)^{2}<0$ if $\beta>\widetilde{\beta}_{p}$.

{\it Case 3:} $p\in(\frac{4}{N},\frac{8}{N})$.
 Let $u_{t}$ be given by \eqref{ut}, then $u_{t}\in \widetilde{S}_{1}$ and  \eqref{eq2.5} \eqref{eq2.6} hold,
Therefore,
\begin{equation} \label{Eut}
\widetilde{E}_{p}^{\beta}(u_{t})=\frac{a}{2}t^{2}+\frac{b}{4}t^{4}-\frac{\beta}{2\|\phi_{p}\|_{L^{2}}^{p}}t^{\frac{Np}{2}}.
\end{equation}
Note now that $\frac{4}{N}<p<\frac{8}{N}$ and $2<\frac{Np}{2}<4$. By Young inequality, we know that
\begin{equation*}
f(t)\triangleq\frac{a}{2}t^{2}+\frac{b}{4}t^{4}\geq\Big(\frac{at^{2}}{2p_{1}}\Big)^{p_{1}}\Big(\frac{bt^{4}}{4q_{1}}\Big)^{q_{1}}
=t^{2+2q_{1}}\Big(\frac{a}{2p_{1}}\Big)^{p_{1}}\Big(\frac{b}{4q_{1}}\Big)^{q_{1}},
\end{equation*}
where $p_{1}+q_{1}$=1 and the equality holds if and only if $\frac{at^{2}}{2p_{1}}=\frac{bt^{4}}{4q_{1}}$. If we let $2q_{1}=\frac{Np}{2}-2$, then $2p_{1}=4-\frac{Np}{2}$ and
\begin{equation}\label{eq2.9}
f(t)\geq t^{\frac{Np}{2}}\Big(\frac{2a}{8-Np}\Big)^{\frac{8-Np}{4}}\Big(\frac{b}{Np-4}\Big)^{\frac{Np-4}{4}}.
\end{equation}
So, in the case of $p\in(\frac{4}{N},\frac{8}{N})$,  we define
\begin{equation}\label{eq2.10}
\widetilde{\beta}_{p}=2\|\phi_{p}\|_{L^{2}}^{p}\Big(\frac{2a}{8-Np}\Big)^{\frac{8-Np}{4}}\Big(\frac{b}{Np-4}\Big)^{\frac{Np-4}{4}}.
\end{equation}
Then, there are  two different situations have to be considered.\\
$\bullet$ If $\beta>\widetilde{\beta}_{p}$, we choose $t_{0}$ such that $\frac{at_{0}^{2}}{2p_{1}}=\frac{bt_{0}^{4}}{4q_{1}}$, then
\begin{equation*}
\widetilde{d}_{\beta}(p)\leq \widetilde{E}_{p}^{\beta}(u_{t_{0}})=\frac{a}{2}t_{0}^{2}+\frac{b}{4}t_{0}^{4}-\frac{\beta}{2\|\phi_{p}\|_{L^{2}}^{p}}t_{0}^{\frac{Np}{2}}
=\frac{t_{0}^{\frac{Np}{2}}}{2\|\phi_{p}\|_{L^{2}}^{p}}(\widetilde{\beta}_{p}-\beta)<0.
\end{equation*}
$\bullet$ If $\beta \leq \widetilde{\beta}_{p}$, it follows from \eqref{eq2.9} and \eqref{eq2.10} that, for any $u\in \widetilde{S}_{1}$,
\begin{equation*}
\aligned
\widetilde{E}_{p}^{\beta}(u)
&=\frac{a}{2}\int_{\mathbb R^{N}}|\nabla u|^{2}dx
+\frac{b}{4}\left(\int_{\mathbb R^{N}}|\nabla u|^{2}dx\right)^{2}
-\frac{\beta}{p+2}\int_{\mathbb R^{N}}|u|^{p+2}dx\\
&\geq\frac{\widetilde{\beta}_{p}-\beta}{2\|\phi_{p}\|_{L^{2}}^{p}}\left(\int_{\mathbb R^{N}}|\nabla u|^{2}dx\right)^{\frac{Np}{4}}\geq0,
\endaligned
\end{equation*}
that is,  $\widetilde{d}_{\beta}(p)\geq0$. Using the fact of \eqref{eqabc}, we see that $\widetilde{d}_{\beta}(p)=0$ for any $\beta \leq \widetilde{\beta}_{p}$.\\
So,  the lemma is proved by combining all the above cases. $\hfill  \square$
\begin{lemma}\label{Lemma12}
If $1\leq N\leq4$ and $p=p^{\ast}=\frac{8}{N}$.  Then, $\widetilde{d}_{\beta}(p^{\ast})=0$ for all $\beta\leq \widetilde{\beta}_{p^{\ast}}$, and $\widetilde{d}_{\beta}(p^{\ast})=-\infty$ for all $\beta>\widetilde{\beta}_{p^{\ast}}$, where $\widetilde{\beta}_{p^{\ast}}$ is defined by \eqref{w11}.
\end{lemma}
Proof: For $1\leq N\leq3$,  this lemma can be proved similarly to that of \cite[Lemma 2.4]{Y1} where $N=4$ is not allowed. In fact, for $1\leq N\leq3$, simply replacing $p$ and $\phi_p$ in \eqref{ut} by $p^*$ and $\phi_{p^*}$ respectively,  we still have \eqref{eq2.5}-\eqref{Eut}, but we note now that $p=p^*=8/N$ and $\frac{Np}{2}=4$, then \eqref{Eut} becomes
\begin{equation*}
\widetilde{E}_{p}^{\beta}(u_{t})=\frac{a}{2}t^{2}+(\frac{b}{4} -\frac{\beta}{2\|\phi_{p}\|_{L^{2}}^{p}})t^4, \text{ for any } t>0,
\end{equation*}
using this fact and the definition of  $\widetilde{\beta}_{p^{\ast}}$ in \eqref{w11}, it is easy to see that the lemma is true for $1\leq N\leq3$.\\
However, if $N=4$, the power $p=p^*=8/N=2$ becomes the critical Sobolev exponent, in this case, although  we still have unique solution $\phi_{p^*}$ (up to translations) for equation \eqref{eq1.9}, but  now $\phi_{p^*}\not\in L^2(\mathbb{R}^4)$ and the above procedures for $1\leq N\leq3$ do not work anymore. So, when $N=4$, we have to redefine $\widetilde{\beta}_{p^{\ast}}$ as in \eqref{w11}.
Then,
for any $\beta\in(0,\widetilde{\beta}_{p^{\ast}}]$ and $p^{\ast}=\frac{8}{N}=2$, by using \eqref{w12} we have
\begin{equation*}
\widetilde{E}_{p^{\ast}}^{\beta}(u)\geq\frac{(\widetilde{\beta}_{p^{\ast}}-\beta)}{4S^{2}}\left(\int_{\mathbb R^{N}}|\nabla u|^{2}dx\right)^{2}\geq0, \text{ for all } u \in S_1,
\end{equation*}
this together with \eqref{eq2.2} imply that $\widetilde{d}_{\beta}(p^{\ast})=0$ for all $\beta \leq \widetilde{\beta}_{p^{\ast}}$.

On the other hand, for any $\beta>\widetilde{\beta}_{p^{\ast}}$, let
\begin{equation*}
U(x)=\frac{2\sqrt{2}}{1+|x|^{2}},\ x\in\mathbb R^{4}.
\end{equation*}
By  \cite[Theorem 1.42]{WM}, we know that $U(x)$ is a minimizer for $S$ and $U(x)$ satisfies
\begin{equation}\label{w13}
\int_{\mathbb R^{4}}|\nabla U(x)|^{2}dx=\int_{\mathbb R^{4}}|U(x)|^{4}dx=S^{2},
\end{equation}
where $S$ is defined by \eqref{w12}. Taking $\eta(x)\in C_{0}^{\infty}(\mathbb R^{4})$ and $0\leq\eta(x)\leq1$ such that $\eta(x)\equiv1$ if $|x|\leq1$, $\eta(x)\equiv0$ if $|x|\geq2$, and $|\nabla\eta(x)|\leq C_{0}$.  Letting
\begin{equation*}
u_{\tau}(x)=A_{\tau}\tau^{2}\eta(x)U(\tau x),
\end{equation*}
where $A_{\tau}$ is chosen so that $\|u_{\tau}\|_{L}^{2}=1$. Then, we have
\begin{equation}\label{w14}
\int_{\mathbb R^{4}}|u_{\tau}(x)|^{2}dx=A_{\tau}^{2}\int_{\mathbb R^{4}}\eta^{2}(\frac{x}{\tau})U^{2}(x)dx=1.
\end{equation}
Since $U(x)\notin L^{2}(\mathbb R^{4})$, there exists a constant $M>0$ such that
\begin{equation}\label{w15}
\int_{\mathbb R^{4}}\eta^{2}(\frac{x}{\tau})U^{2}(x)dx\geq\int_{|x|\leq\tau}U^{2}(x)dx\rightarrow+\infty\ \mbox{as}\ \tau\rightarrow+\infty
\end{equation}
and
\begin{equation}\label{w16}
\int_{\mathbb R^{4}}\eta^{2}(\frac{x}{\tau})U^{2}(x)dx\leq\int_{|x|\leq2\tau}U^{2}(x)dx\leq M\ln2\tau\ \mbox{as}\ \tau\ \mbox{large enough}.
\end{equation}
By \eqref{w14}-\eqref{w16}, we know that
\begin{equation}\label{w17}
A_{\tau}^{2}\rightarrow0\ \mbox{as}\ \tau\rightarrow+\infty\quad\mbox{and}\quad A_{\tau}^{2}\ln2\tau\geq \frac{1}{M}\ \mbox{as}\ \tau\ \mbox{large enough}.
\end{equation}
Using \eqref{w13} and \eqref{w17}, we have
\begin{equation}\label{w18}
\int_{\mathbb R^{4}}u_{\tau}^{4}(x)dx=A_{\tau}^{4}\tau^{4}\int_{\mathbb R^{4}}\eta^{4}(\frac{x}{\tau})U^{4}(x)dx\geq
A_{\tau}^{4}\tau^{4}\int_{|x|\leq\tau}U^{4}(x)dx,
\end{equation}
\begin{equation}\label{w19}
\int_{\mathbb R^{4}}|\nabla u_{\tau}(x)|^{2}dx=A_{\tau}^{2}\tau^{4}\int_{\mathbb R^{4}}|\nabla\eta(x)U(\tau x)+\tau\eta(x)\nabla U(\tau x)|^{2}dx\\
\leq A_{\tau}^{2}\tau^{2}S^{2}+O(A_{\tau}^{2}\tau).
\end{equation}
Since $\beta>\widetilde{\beta}_{p^{\ast}}$,  it follows from \eqref{w13} and \eqref{w17}-\eqref{w19} that
\begin{equation}\label{w41}
\widetilde{d}_{\beta}(p^{\ast})\leq\widetilde{E}_{p^{\ast}}^{\beta}(u_{\tau})\leq\frac{a}{2}A_{\tau}^{2}\tau^{2}S^{2}
+\frac{A_{\tau}^{4}\tau^{4}}{4}\left(bS^{4}-\beta\int_{|x|\leq\tau}U^{4}(x)dx\right)
+O(A_{\tau}^{4}\tau^{3}).
\end{equation}
Then, by \eqref{w13} and \eqref{w17}, we have
\begin{equation}\label{w42}
A_{\tau}^{2}\tau=A_{\tau}^{2}\ln2\tau\frac{\tau}{\ln2\tau}\rightarrow+\infty\quad\mbox{and}\quad
bS^{4}-\beta\int_{|x|\leq\tau}U^{4}(x)dx\rightarrow S^{2}(bS^{2}-\beta)<0\ \mbox{as}\ \tau\rightarrow+\infty.
\end{equation}
Hence,  by \eqref{w41}, \eqref{w42} and let $\tau\rightarrow+\infty$, we see that $\widetilde{d}_{\beta}(p^{\ast})=-\infty$ for all $\beta>\widetilde{\beta}_{p^{\ast}}$.
$\hfill \square$


\begin{lemma}\label{le3}
Let $u\in H^{1}(\mathbb R^{N})$ $(N \geq 1)$, then there exists a nonnegative, non-increasing function $u^{\ast}\in H_{r}^{1}(\mathbb R^{N})$
such that
\begin{equation*}
\int_{\mathbb R^{N}}|u^{\ast}|^{p}dx=\int_{\mathbb R^{N}}|u|^{p}dx
\end{equation*}
for all $1\leq p<\infty$ and
\begin{equation*}
\int_{\mathbb R^{N}}|\nabla u^{\ast}|^{2}dx\leq\int_{\mathbb R^{N}}|\nabla u|^{2}dx.
\end{equation*}
\end{lemma}
Proof. The proof of the lemma  can be found in  \cite[ appendix A.III]{BL}. $\hfill \square$
\begin{lemma}\label{le2}
\cite[Proposition 1.7.1]{CH}
 Let $\{u_{n}\}\subset H_{r}^{1}(\mathbb R^{N})$ be a bounded sequence. If $N\geq2$ or if $u_{n}(x)$ is a non-increasing  function of $|x|$ for
every $n\geq0$, then there exist a subsequence $\{u_{n_{k}}\}_{k\geq0}$ and $u\in H_{r}^{1}(\mathbb R^{N})$ such that $u_{n_{k}}\rightarrow u$ as
$k\rightarrow\infty$ in $L^{q}(\mathbb R^{N})$ for $q\in(2,\frac{2N}{N-2})$ if $N\geq3$, or $q\in(2,+\infty)$ if $N=1,2$.
\end{lemma}
\begin{lemma}\label{le6}
If $V(x)\equiv0$, $p\in(0,\frac{8}{N})$, $1\leq N\leq4$ and $\widetilde{d}_{\beta}(p)<0$. Then,
\eqref{eq1.3} has at least a nonnegative  minimizer.
\end{lemma}
Proof. Let $\{u_{n}\}\subset \widetilde{S}_{1}$ be a minimizing sequence of $\widetilde{d}_{\beta}(p)$,  then it is easy to know that $\{u_{n}\}$ is bounded in $H^{1}(\mathbb R^{N})$ by using \eqref{a4} and $\widetilde{d}_{\beta}(p)<0$.
By Lemma \ref{le3}, we know that there exists $\{u_{n}^{\ast}\}\subset H_{r}^{1}(\mathbb R^{N})$ which are nonnegative, non-increasing function and $\{u_{n}^{\ast}\}\subset \widetilde{S}_{1}$
is also a minimizing sequence for $\widetilde{d}_{\beta}(p)$. Moreover, $\{u_{n}^{\ast}\}$ is still bounded in $H_{r}^{1}(\mathbb R^{N})$. By Lemma \ref{le2}, there exist a subsequence, still denoted by $\{u_{n}^{\ast}\}$, and some $u_{0}\in H_{r}^{1}(\mathbb R^{N})$ such that
\begin{equation}\label{eq3.1}
u_{n}^{\ast}\overset{n}{\rightharpoonup} u_{0} \mbox{ weakly in } \ H_{r}^{1}(\mathbb R^{N}), \text{ and }
u_{n}^{\ast}\overset{n}{\rightarrow} u_{0}\quad\mbox{strongly in}\  L^{q}(\mathbb R^{N})
\end{equation}
for $q\in(2,\frac{2N}{N-2})$ if $N\geq3$, or $q\in(2,+\infty)$ if $N=1,2$.

We claim that $u_0 \not\equiv 0$. Otherwise, \eqref{eq3.1} implies that
\begin{equation}\label{eq3.2}
\int_{\mathbb R^{N}}|u_{n}^{\ast}|^{p+2}dx \rightarrow 0 \text{ as } n\rightarrow +\infty.
\end{equation}
Then,  by the definition of \eqref{Ebeta} we know that
\begin{equation*}
0\leq\lim_{k\rightarrow\infty}\left[\frac{a}{2}\int_{\mathbb R^{N}}|\nabla u_{n_{k}}^{\ast}|^{2}dx+
\frac{b}{4}\left(\int_{\mathbb R^{N}}|\nabla u_{n_{k}}^{\ast}|^{2}dx\right)^{2}\right]=\widetilde{d}_{\beta}(p)<0,
\end{equation*}
a contradiction. Hence,
$u_{0}\not\equiv 0$ 
and
\begin{equation}\label{eq3.3}
\widetilde{E}_{p}^{\beta}(u_{0})\leq\lim_{n\rightarrow\infty}\widetilde{E}_{p}^{\beta}(u_{n}^{\ast})=\widetilde{d}_{\beta}(p)<0.
\end{equation}
Let $\gamma=\|u_{0}\|_{L^{2}}^{2} $, then $\gamma \in (0,1]$ and  $u_{\gamma}(x):=u_{0}(\gamma^{\frac{1}{N}}x) \in \widetilde{S}_{1}$. Hence, it follows from  \eqref{eq3.3} that
\begin{equation} \label{gamma1}
\aligned
\widetilde{d}_{\beta}(p)\leq \widetilde{E}_{p}^{\beta}(u_{\gamma})
&=\frac{a\gamma^{\frac{2}{N}-1}}{2}\int_{\mathbb R^{N}}|\nabla u_{0}|^{2}dx
+\frac{b\gamma^{\frac{4}{N}-2}}{4}\left(\int_{\mathbb R^{N}}|\nabla u_{0}|^{2}dx\right)^{2}
-\frac{\beta}{\gamma (p+2)}\int_{\mathbb R^{N}}|u_{0}|^{p+2}dx\\
&=\frac{1}{\gamma}\left[\frac{a\gamma^{\frac{2}{N}}}{2}\int_{\mathbb R^{N}}|\nabla u_{0}|^{2}dx
+\frac{b\gamma^{\frac{4}{N}-1}}{4}\left(\int_{\mathbb R^{N}}|\nabla u_{0}|^{2}dx\right)^{2}
-\frac{\beta}{p+2}\int_{\mathbb R^{N}}|u_{0}|^{p+2}dx\right]\\
&\leq\frac{1}{\gamma}\widetilde{E}_{p}^{\beta}(u_{0})\leq\frac{1}{\gamma}\widetilde{d}_{\beta}(p),
\endaligned
\end{equation}
this implies that $\frac{1}{\gamma}\leq1$ since $\widetilde{d}_{\beta}(p)<0$. Then, $\gamma\geq1$ and $\|u_{0}\|_{L^{2}}=1$. Hence, \eqref{eq3.3} implies that
$u_{0}$ is a minimizer of $\widetilde{d}_{\beta}(p)$. Moreover,  by the definition of  $\widetilde{d}_{\beta}(p)$, we know that $|u_{0}|$ is also a minimizer,  so we may assume that $\widetilde{d}_{\beta}(p)$ has a nonnegative  minimizer.
$  \hfill\square$

Now, we are ready to prove our Theorem \ref{th5}\\
{\bf Proof of Theorem \ref{th5} :}
$(\textup{i})$  When $p\in (0, 8/N)$, it follows from Lemma \ref{lea} that $\widetilde{d}_{\beta}(p)=0$ for all $ \beta \leq \widetilde{\beta}_{p}$, which then shows that $\widetilde{d}_{\beta}(p)$ has no any minimizer for all $ \beta < \widetilde{\beta}_{p}$. Otherwise, if there exists $ \beta < \widetilde{\beta}_{p}$ such that $\widetilde{d}_{\beta}(p)$ has a minimizer $u\in S_1$, that is,
$\widetilde{E}_{p}^{\beta}(u)=\widetilde{d}_{\beta}(p)=0$, and then, by $ \beta < \widetilde{\beta}_{p}$ and the definition of \eqref{eq1.7} we see that
\begin{equation*}
0=\widetilde{d}_{\widetilde{\beta}_{p}}(p)\leq \widetilde{E}_{p}^{\widetilde{\beta}_{p}}(u)<\widetilde{E}_{p}^{\beta}(u)=0,
\end{equation*}
which is impossible. Particularly, if $p\in (0, 4/N]$, we claim that $\widetilde{d}_{\beta}(p)$ has no minimizer even for $ \beta = \widetilde{\beta}_{p}$. In fact, by the definition \eqref{a2} we know that $ \widetilde{\beta}_{p}=0$ for $p\in (0, 4/N)$ and $ \widetilde{\beta}_{p}>0$ for $p=4/N$. If $\beta =\widetilde{\beta}_{p}=0$ and  there is a minimizer $u_0$ for $\widetilde{d}_{\beta}(p)$,
we then have $u_0 \equiv 0$ by  using \eqref{a4} and the fact  that $\widetilde{d}_{\beta}(p) = \widetilde{E}_{p}^{\beta}(u_{0})=0$, which leads to a contradiction  since $u_{0}\in \widetilde{S}_{1}$.
On the other hand,
 for $p=\frac{4}{N}$, we know also that $\widetilde{d}_{\beta}(p)=0$  for $\beta=\widetilde{\beta}_{p}$, 
 if there exists $u_{0}\in \widetilde{S}_{1}$ such that $\widetilde{E}_{p}^{\beta}(u_{0})=\widetilde{d}_{\beta}(p)=0$. Then,  using  \eqref{eq1.8}  and the value of  $ \widetilde{\beta}_{p}$ for $p=\frac{4}{N}$ in \eqref{a2},    we have
\begin{equation*}
\frac{a}{2}\int_{\mathbb R^{N}}|\nabla u_{0}|^{2}dx+\frac{b}{4}\left(\int_{\mathbb R^{N}}|\nabla u_{0}|^{2}dx\right)^{2}
=\frac{\beta}{2+\frac{4}{N}}\int_{\mathbb R^{N}}|u_{0}|^{2+\frac{4}{N}}dx\leq\frac{a}{2}\int_{\mathbb R^{N}}|\nabla u_{0}|^{2}dx,
\end{equation*}
this implies that   $u_{0}=0$, which contradicts $u_{0}\in \widetilde{S}_{1}$ .

Now, we come to prove the existence.
When $p\in(0,\frac{8}{N})$ and $\beta>\widetilde{\beta}_{p}$,  as a straightforward consequent of  Lemmas \ref{lea} and \ref{le6} we know  that $\widetilde{d}_{\beta}(p)$ has a nonnegative minimizer.
When $\beta=\widetilde{\beta}_{p}$, we know that $\widetilde{d}_{\widetilde{\beta}_{p}}(p) =0$ by Lemma \ref{lea}.
In what follows, we show that, if $p\in(\frac{4}{N},\frac{8}{N})$,  
 $\widetilde{d}_{\widetilde{\beta}_{p}}(p)=0$ has also a minimizer.

Let $\beta_{n}=\widetilde{\beta}_{p}+\frac{1}{n}$ with $1\leq n\in \mathbb{Z}$,  then, for each $\beta_{n}$,  Lemmas \ref{lea} and \ref{le6} tell us that  there exists $u_{n}\in \widetilde{S}_{1}$ such that
\begin{equation}\label{eq3.4}
\widetilde{d}_{\beta_{n}}(p)=\widetilde{E}_{p}^{\beta_{n}}(u_{n})<0,
\end{equation}
and $\{u_{n}\}$ is bounded in $H^{1}(\mathbb R^{N})$. It follows from Lemmas \ref{le3} and \ref{le2}  that
there exists $u_{n}^{\ast} \in H_{r}^{1}(\mathbb R^{N})\cap   \widetilde{S}_{1}$ which
is also a minimizer of $\widetilde{d}_{\beta_{n}}(p)$  and $\{u_{n}^{\ast}\}\subset H_{r}^{1}(\mathbb R^{N})$ is bounded. Hence,  there is $ u_{0}\in H_{r}^{1}(\mathbb R^{N})$ such that
\begin{equation}\label{eq3.5}
u_{n}^{\ast}\overset{n}{\rightharpoonup} u_{0} \mbox{ weakly in } \ H_{r}^{1}(\mathbb R^{N}), \text{ and }
u_{n}^{\ast}\overset{n}{\rightarrow} u_{0}\quad\mbox{strongly in}\quad L^{q}(\mathbb R^{N}),
\end{equation}
for $q\in(2,\frac{2N}{N-2})$ if $N\geq3$, or $q\in(2,+\infty)$ if $N=1,2$.
Since $\beta_{n}=\widetilde{\beta}_{p}+\frac{1}{n}$, by the definition of $\widetilde{E}_{p}^{\beta_{n}}(u_{n}^{\ast})$ and \eqref{eq3.4} we have
\begin{equation}\label{b9}
\widetilde{d}_{\widetilde{\beta}_{p}}(p)-\frac{1}{n(p+2)}\int_{\mathbb R^{N}}|u_{n}^{\ast}|^{p+2}dx
\leq \widetilde{E}_{p}^{\beta_{n}}(u_{n}^{\ast})= \widetilde{d}_{\beta_{n}}(p)<0.
\end{equation}
Note that $\widetilde{d}_{\widetilde{\beta}_{p}}(p)=0$ by Lemma \ref{lea}, we then follows from \eqref{b9} that
\begin{equation}\label{c1}
\widetilde{d}_{\beta_{n}}(p)\rightarrow0\quad\mbox{as}\ n\rightarrow+\infty.
\end{equation}
For $u_0$ given in \eqref{eq3.5}, we claim that $u_0\not \equiv 0$.
Otherwise, if
$\int_{\mathbb R^{N}}|u_{n}^{\ast}|^{p+2}dx\rightarrow0$ as $n\rightarrow\infty$, then
by \eqref{c1},
we have
$\frac{a}{2}\int_{\mathbb R^{N}}|\nabla u_{n }^{\ast}|^{2}dx
+\frac{b}{4}\left(\int_{\mathbb R^{N}}|\nabla u_{n }^{\ast}|^{2}dx\right)^{2}\rightarrow0$ as $n\rightarrow\infty$, that is, $\int_{\mathbb R^{N}}|\nabla u_{n }^{\ast}|^{2}dx \rightarrow0$ as $n\rightarrow\infty$.
Therefore, using \eqref{a4} for $u_n^*$ and $\beta_n$, by $1<Np/4<2$ we know that
$\widetilde{d}_{\beta_{n}}(p)=\widetilde{E}_{p}^{\beta_{n}}(u_{n}^{\ast})>0$ for $n$ large enough, this however contradicts \eqref{b9}. So, $u_{0}\neq0$. Moreover, it follows from \eqref{eq3.5} that $0<\|u_{0}\|_{L^{2}}\leq1$ and
\begin{equation}\label{eq3.8}
\widetilde{E}_{p}^{\widetilde{\beta}_{p}}(u_{0})\leq\lim_{n\rightarrow\infty}\widetilde{E}_{p}^{\beta_{n}}(u_{n}^{\ast})=0.
\end{equation}
Let $\gamma=\|u_{0}\|_{L^{2}}^{2}$, then $\gamma \in (0,1]$ and $u_{\gamma}(x):=u_{0}(\gamma^{\frac{1}{N}}x) \in \widetilde{S}_{1}$. If $\gamma \not =1$, i.e., $0<\gamma<1$, similar to the derivation of \eqref{gamma1} we have
\begin{equation*}
\aligned
\widetilde{d}_{\widetilde{\beta}_{p}}(p)\leq \widetilde{E}_{p}^{\widetilde{\beta}_{p}}(u_{\gamma})
&=\frac{1}{\gamma}\Big[\frac{a\gamma^{\frac{2}{N}}}{2}\int_{\mathbb R^{N}}|\nabla u_{0}|^{2}dx
+\frac{b\gamma^{\frac{4}{N}-1}}{4}\left(\int_{\mathbb R^{N}}|\nabla u_{0}|^{2}dx\right)^{2}
-\frac{\widetilde{\beta}_{p}}{p+2}\int_{\mathbb R^{N}}|u_{0}|^{p+2}dx\Big]\\
&<\frac{1}{\gamma}\widetilde{E}_{p}^{\widetilde{\beta}_{p}}(u_{0})\leq0, \text{ by } \gamma<1 \text{ and  \eqref{eq3.8}},
\endaligned
\end{equation*}
that is, $\widetilde{d}_{\widetilde{\beta}_{p}}(p)<0$. However, $\widetilde{d}_{\widetilde{\beta}_{p}}(p)=0$ by Lemma \ref{lea}. So,  $\|u_{0}\|_{L^{2}}=\gamma=1$ and then, using \eqref{eq3.8}, we know that $u_{0}$ is a minimizer of $\widetilde{d}_{\widetilde{\beta}_{p}}(p)(=0)$.

If $ {u} $ is a minimizer of $\widetilde{d}_{\beta}(p)$, it is well-known that  there is a lagrange multiplier $\lambda \in\mathbb R$ such that
\begin{equation*} 
-\left(a+b\int_{\mathbb R^{N}}|\nabla  {u} |^{2}dx\right)\triangle  {u} =\beta| {u} |^{p} {u} +\lambda u .
\end{equation*}
Since $\|u\|_{L^2}=1$, multiplying $u$ both sides in the above  equation and integrating ,  we have
\begin{equation}\label{ww}
\lambda =a\int_{\mathbb R^{N}}|\nabla  {u} |^{2}dx
+b\left(\int_{\mathbb R^{N}}|\nabla  {u} |^{2}dx\right)^{2}
-\beta\int_{\mathbb R^{N}}| {u} |^{p+2}dx.
\end{equation}
Moreover, $u$ satisfies the following Pohozaev identity \cite{LY}
\begin{equation*}
\frac{a(N-2)}{2}\int_{\mathbb R^{N}}|\nabla u|^{2}dx
+\frac{b(N-2)}{2}\left(\int_{\mathbb R^{N}}|\nabla u|^{2}dx\right)^{2}
-\frac{\beta N}{p+2}\int_{\mathbb R^{N}}|u|^{p+2}dx=\frac{N}{2}\lambda,
\end{equation*}
which together with \eqref{ww} imply that
\begin{equation*}
\lambda=\frac{(N-2)\beta p-4\beta}{2(p+2)}\int_{\mathbb R^{N}}|u|^{p+2}dx<0.
\end{equation*}
This shows that, for some $\lambda<0$,  \eqref{eq1.1} has a nonnegative solution $u\geq 0$ with $\|u\|_{L^2}=1$  and we know also that  $u > 0$ by the strong maximum principle.\\
$(\textup{ii})$ For $p=\frac{8}{N}$, by Lemma \ref{Lemma12} it is clear that $\widetilde{d}_{\beta}(p^{\ast})$ has no minimizer for any $\beta>\widetilde{\beta}_{p^{\ast}}$. Suppose that there exists some $\beta\leq \widetilde{\beta}_{p^{\ast}}$ such that $\widetilde{d}_{\beta}(p^{\ast})$ has a  minimizer $u_{0}\in \widetilde{S}_{1}$, then by $\widetilde{d}_{\beta}(p^{\ast})=0$ (Lemma \ref{Lemma12}) and \eqref{w12} if $N=4$, or \eqref{eq1.8} if $N\leq 3$, we have
\begin{equation*}
\frac{a}{2}\int_{\mathbb R^{N}}|\nabla u_{0}|^{2}dx+\frac{b}{4}\left(\int_{\mathbb R^{N}}|\nabla u_{0}|^{2}dx\right)^{2}
=\frac{\beta}{2+\frac{8}{N}}\int_{\mathbb R^{N}}|u_{0}|^{2+\frac{8}{N}}dx\leq\frac{b}{4}\left(\int_{\mathbb R^{N}}|\nabla u_{0}|^{2}dx\right)^{2},
\end{equation*}
that is, $u_{0}\equiv0$, which is impossible. Hence, $\widetilde{d}_{\beta}(p^{\ast})$ has no minimizer for all $\beta\leq \widetilde{\beta}_{p^{\ast}}$, either.$\hfill \square$

In the end of this section, we give some estimates on $\widetilde{d}_{\beta}(p)$, which are required in section 3.
\begin{lemma}\label{Le3.1}
For $p\in (0,p^*)$ and $\beta_{p}$ in \eqref{eq1.10},
let $\widetilde{d}_{\beta}(p)$ be defined in \eqref{eq1.7}. Then, for any fixed $\beta>\beta_{p^{\ast}}$ and $1\leq N\leq3$, we have
\begin{equation*}
\widetilde{d}_{\beta}(p)=-\frac{b(p^{\ast}-p)}{4p}\left(\frac{\beta p}{\beta_{p}p^{\ast}}\right)^{\frac{p^{\ast}}{p^{\ast}-p}}(1+o(1))\quad and\quad \widetilde{d}_{\beta}(p)\rightarrow-\infty\quad \mbox{as}\ p\nearrow p^{\ast}.
\end{equation*}
\end{lemma}
Proof. Let $ u\in \widetilde{S}_{1}$,  then , it follows from \eqref{eq1.10} and \eqref{eq1.8}  that
\begin{equation}\label{b1}
\aligned
\widetilde{E}_{p}^{\beta}(u)
&=\frac{a}{2}\int_{\mathbb R^{N}}|\nabla u|^{2}dx
+\frac{b}{4}\left(\int_{\mathbb R^{N}}|\nabla u|^{2}dx\right)^{2}
-\frac{\beta}{p+2}\int_{\mathbb R^{N}}|u|^{p+2}dx\\
&\geq\frac{b}{4}\left(\int_{\mathbb R^{N}}|\nabla u|^{2}dx\right)^{2}
-\frac{\beta}{2\|\phi_{p}\|_{L^{2}}^{p}}\left(\int_{\mathbb R^{N}}|\nabla u|^{2}dx\right)^{\frac{Np}{4}}\\
&=\frac{b}{4}\left(\int_{\mathbb R^{N}}|\nabla u|^{2}dx\right)^{2}
-\frac{b\beta}{4\beta_{p}}\left(\int_{\mathbb R^{N}}|\nabla u|^{2}dx\right)^{\frac{Np}{4}}.
\endaligned
\end{equation}
Denoting
\begin{equation}\label{c3}
\left(\int_{\mathbb R^{N}}|\nabla u|^{2}dx\right)^{2}=r\quad\mbox{and}\quad h(r)\triangleq\frac{b}{4}r-\frac{b\beta}{4\beta_{p}}r^{\frac{p}{p^{\ast}}}(r>0),
\end{equation}
 by simple computation, we know that the function $h$ has a unique minimum at $r:=r_p$ with
\begin{equation}\label{eq3.7}
r_{p}\triangleq\left(\frac{\beta p}{\beta_{p}p^{\ast}}\right)^{\frac{p^{\ast}}{p^{\ast}-p}}.
\end{equation}
That is,
\begin{equation}\label{b2}
h(r)\geq h(r_{p})=-\frac{b(p^{\ast}-p)}{4p}\left(\frac{\beta p}{\beta_{p}p^{\ast}}\right)^{\frac{p^{\ast}}{p^{\ast}-p}}.
\end{equation}
Then,  \eqref{b1} and \eqref{b2} shows that
\begin{equation*}
\widetilde{d}_{\beta}(p)\geq-\frac{b(p^{\ast}-p)}{4p}\left(\frac{\beta p}{\beta_{p}p^{\ast}}\right)^{\frac{p^{\ast}}{p^{\ast}-p}},
\end{equation*}
and we have a lower bound for $\widetilde{d}_{\beta}(p)$. Now, we come to estimate the upper bound for $\widetilde{d}_{\beta}(p)$.

Let $u_{t}=\frac{t^{\frac{N}{2}}\phi_{p}(tx)}{\|\phi_{p}\|_{L^{2}}}$ for $t>0$, then $\int_{\mathbb R^{N}}u_{t}^{2}dx=1$. Similar to \eqref{eq2.5} and \eqref{eq2.6}, we have
\begin{equation*}
\aligned
\widetilde{E}_{p}^{\beta}(u_{t})&=\frac{a}{2}t^{2}+\frac{b}{4}t^{4}
-\frac{b\beta}{4\beta_{p}}\left(t^{4}\right)^{\frac{p}{p^{\ast}}}.
\endaligned
\end{equation*}
Taking $t^{4}=\left(\frac{\beta p}{\beta_{p}p^{\ast}}\right)^{\frac{p^{\ast}}{p^{\ast}-p}}$, then
\begin{equation}\label{b6}
\widetilde{E}_{p}^{\beta}(u_{t})
=\frac{a}{2}\left(\frac{\beta p}{\beta_{p}p^{\ast}}\right)^{\frac{p^{\ast}}{2(p^{\ast}-p)}}
-\frac{b(p^{\ast}-p)}{4p}\left(\frac{\beta p}{\beta_{p}p^{\ast}}\right)^{\frac{p^{\ast}}{p^{\ast}-p}}.
\end{equation}
By \eqref{eq1.13}, we know that
\begin{equation*}
\frac{\frac{a}{2}\left(\frac{\beta p}{\beta_{p}p^{\ast}}\right)^{\frac{p^{\ast}}{2(p^{\ast}-p)}}}
{\frac{b(p^{\ast}-p)}{4p}\left(\frac{\beta p}{\beta_{p}p^{\ast}}\right)^{\frac{p^{\ast}}{p^{\ast}-p}}}\rightarrow0\quad\mbox{as}\ p\nearrow p^{\ast}.
\end{equation*}
Then, using \eqref{b6} we see  that
\begin{equation*}
\widetilde{d}_{\beta}(p)
\leq\widetilde{E}_{p}^{\beta}(u_{t})
=-\frac{b(p^{\ast}-p)}{4p}\left(\frac{\beta p}{\beta_{p}p^{\ast}}\right)^{\frac{p^{\ast}}{p^{\ast}-p}}(1+o(1)),
\end{equation*}
 So, we finish the proof of the  lemma.$\quad\quad\quad\quad\quad\quad\quad\quad\quad\quad\quad\quad\quad\quad\quad\quad\quad\square$
\section{Case of $V(x)\not\equiv0$}
In this section, we come to prove Theorems \ref{th1.1}-\ref{th2}. For this purpose, we first recall an embedding theorem which can be found in \cite[Theorem XIII.67]{RS} or \cite[Theorem 2.1]{BW}, etc.
\begin{lemma}\label{Lm1}
For any $N\geq1$, let $V(x)$ satisfy the condition \eqref{a3}, then the embedding from $\mathcal{H}$ into $L^{q}(\mathbb R^{N})(2\leq q<2^{\ast})$ is compact.
$\quad\quad\quad\quad\quad\quad\quad\quad\quad\quad\quad\quad
\quad\quad\quad\quad\quad\quad\quad\quad\quad\quad\quad\quad\quad\quad\quad\quad\quad\quad\square$
\end{lemma}
{\bf Proof of Theorem \ref{th1.1} :} For any fixed $\beta>0$, $0<p<p^{\ast}$ and $u\in S_{1},$ it follows from \eqref{eq1.10} and \eqref{eq1.8} that
\begin{equation}\label{o1}
\aligned
E_{p}^{\beta}(u)
&=\frac{a}{2}\int_{\mathbb R^{N}}|\nabla u|^{2}dx+\frac{b}{4}\left(\int_{\mathbb R^{N}}|\nabla u|^{2}dx\right)^{2}
+\frac{1}{2}\int_{\mathbb R^{N}}V(x)u^{2}dx-\frac{\beta}{p+2}\int_{\mathbb R^{N}}|u|^{p+2}dx\\
&\geq\frac{a}{2}\int_{\mathbb R^{N}}|\nabla u|^{2}dx
+\frac{b}{4}\left(\int_{\mathbb R^{N}}|\nabla u|^{2}dx\right)^{2}+\frac{1}{2}\int_{\mathbb R^{N}}V(x)u^{2}dx
-\frac{\beta}{2\|\phi_{p}\|_{L^{2}}^{p}}\left(\int_{\mathbb R^{N}}|\nabla u|^{2}dx\right)^{\frac{Np}{4}}.
\endaligned
\end{equation}
Since $ 0<\frac{Np}{4}<2$,  using \eqref{o1} it is easy to see that $d_{\beta}(p)>-\infty$. Let $\{u_{n}\}$ be a minimizing sequence of $d_{\beta}(p)$, then it is not difficult to know that
$\{u_{n}\}$ is bounded in $\mathcal{H}$. Then,
by Lemma \ref{Lm1}, for some $u\in S_{1}$,  we may assume that,  passing to a subsequence if necessary,
\begin{equation*}
u_{n}\rightharpoonup u  \mbox{ weakly in }  \mathcal{H},\
u_{n}\rightarrow u  \mbox{ strongly in }  L^{q}(\mathbb R^{N}) \text{ with } q \in [ 2, 2^{\ast})
\end{equation*}
as $n\rightarrow\infty$, and
\begin{equation*}
\left(\int_{\mathbb R^{N}}|\nabla u|^{2}dx\right)^{2}\leq\liminf_{n\rightarrow\infty}\left(\int_{\mathbb R^{N}}|\nabla u_{n}|^{2}dx\right)^{2} .
\end{equation*}
Hence
\begin{equation*}
d_{\beta}(p)=\liminf_{n\rightarrow\infty} E_{p}^{\beta}(u_{n})\geq E_{p}^{\beta}(u)\geq d_{\beta}(p).
\end{equation*}
Therefore, $u$ is a minimizer of $d_{\beta}(p)$ for all $\beta>0$.$\quad\quad\quad\quad\quad
\quad\quad\quad\quad\quad\quad\quad\quad\quad\quad\quad\quad\quad\quad\quad\quad\quad\quad\quad\square$\\
{\bf Proof of Theorem \ref{th1} :} (\textup{i})
Taking $p=p^*=8/N$ and $u\in S_{1}$ in \eqref{o1},  by the definition of $\beta_{p^*}$ and \eqref{eq1.8},  it is easy to see that there exists some constant $c>0$ such that
\begin{align} \label{p*}
E_{p^*}^{\beta}(u)
&\geq\frac{a}{2}\int_{\mathbb R^{N}}|\nabla u|^{2}dx
+\frac{1}{2}\int_{\mathbb R^{N}}V(x)u^{2}dx
+\frac{\beta_{p^*}-\beta}{2\|\phi_{p^*}\|_{L^{2}}^{p^*}}\left(\int_{\mathbb R^{N}}|\nabla u|^{2}dx\right)^2 \nonumber\\
&\geq\frac{a}{2}\int_{\mathbb R^{N}}|\nabla u|^{2}dx
+\frac{1}{2}\int_{\mathbb R^{N}}V(x)u^{2}dx \geq c\|u\|^2_{L^2} =c>0,
\end{align}
since  $\beta \in (0, \beta_{p^*}] $  and the embedding lemma \ref{Lm1}. This shows that $d_{\beta}(p^*)>0$. Now, we come to prove that $d_{\beta}(p^*)$ can be attained for all  $\beta \in (0, \beta_{p^*}] $. In fact, let $\{u_n\}$ be a minimizing sequence for $d_{\beta}(p^*)$, then, by \eqref{p*} we know that $\{u_n\}$ is bounded in $\mathcal{H}$ and $\{u_{n}\}$ converges weakly in $\mathcal{H}$ to some $u\in S_{1}$ as in the proof of Theorem \ref{th1.1}. Therefore,
\begin{equation*}
d_{\beta}(p^*)=\liminf_{n\rightarrow\infty} E_{p^*}^{\beta}(u_{n}) \geq E_{p^*}^{\beta}(u) \geq d_{\beta}(p^*),
\end{equation*}
That is, $u$ is a minimizer of $d_{\beta}(p^*)$ for all $\beta \in (0,\beta_{p^*}]$.

(\textup{ii}) Let $\varphi(x)\in C_{0}^{\infty}(\mathbb R^{N})$ and $0\leq\varphi(x)\leq1$ such that $\varphi(x)\equiv1$ if $|x|\leq1$, $\varphi(x)\equiv0$ if $|x|\geq2$, and $|\nabla\varphi(x)|\leq 2$. For any $x_{0}\in\mathbb{R}^{N}$, we set
\begin{equation*}
u_{\tau}=\frac{A_{\tau}\tau^{\frac{N}{2}}}{\|\phi_{p^{\ast}}\|_{L^{2}}}\varphi(x-x_{0})\phi_{p^{\ast}}(\tau(x-x_{0})),
\end{equation*}
where $A_{\tau}>0$ is chosen so that $\|u_{\tau}\|_{L^{2}}=1$.
By the exponential decay of $\phi_{p^{\ast}}$ (see, e.g. \cite{GNN}), we know that
\begin{equation*}
\frac{1}{A_{\tau}^{2}}=
1+\frac{1}{\|\phi_{p^{\ast}}\|_{L^{2}}^{2}}\int_{\mathbb{R}^{N}}\left(\varphi^{2}(\frac{x}{\tau})-1\right)\phi_{p^{\ast}}^{2}(x)dx\rightarrow1\quad\mbox{as}\ \tau\rightarrow\infty.
\end{equation*}
Then,  $A_{\tau}\geq1$ and $\lim\limits_{\tau\rightarrow+\infty}A_{\tau}=1$.
Moreover, it follows from \eqref{eq1.10} and \eqref{eq1.11} that
\begin{equation}\label{eq2.1}
\aligned
\int_{\mathbb R^{N}}|\nabla u_{\tau}|^{2}dx&=\frac{A_{\tau}^{2}\tau^{N}}{\|\phi_{p^{\ast}}\|_{L^{2}}^{2}}
\int_{\mathbb{R}^{N}}|\nabla\varphi(x-x_{0})\phi_{p^{\ast}}(\tau(x-x_{0}))+\tau\varphi(x-x_{0})\nabla\phi_{p^{\ast}}(\tau(x-x_{0}))|^{2}dx\\
&=A_{\tau}^{2}\tau^{2}+O(\tau^{-\infty}),
\endaligned
\end{equation}
\begin{equation}\label{a5}
\aligned
\int_{\mathbb R^{N}}|u_{\tau}|^{p^{\ast}+2}dx&=\frac{A_{\tau}^{p^{\ast}+2}\tau^{N+4}}{\|\phi_{p^{\ast}}\|_{L^{2}}^{p^{\ast}+2}}
\int_{\mathbb R^{N}}\varphi^{p^{\ast}+2}(x-x_{0})\phi_{p^{\ast}}^{p^{\ast}+2}(\tau(x-x_{0}))dx\\
&=\frac{bA_{\tau}^{p^{\ast}+2}\tau^{4}(p^{\ast}+2)}{4\beta_{p^{\ast}}}+O(\tau^{-\infty}),
\endaligned
\end{equation}
\begin{equation}\label{a6}
\aligned
\int_{\mathbb R^{N}}V(x)u_{\tau}^{2}dx&=\frac{A_{\tau}^{2}\tau^{N}}{\|\phi_{p^{\ast}}\|_{L^{2}}^{2}}
\int_{\mathbb R^{N}}V(x)\varphi^{2}(x-x_{0})\phi_{p^{\ast}}^{2}(\tau(x-x_{0}))dx\\
&=\frac{A_{\tau}^{2}}{\|\phi_{p^{\ast}}\|_{L^{2}}^{2}}
\int_{\mathbb R^{N}}V(\frac{x}{\tau}+x_{0})\varphi^{2}(\frac{x}{\tau})\phi_{p^{\ast}}^{2}(x)dx.\\
&=V(x_{0})+o(1),
\endaligned
\end{equation}
where, and in what follows,  the notation $O(\tau^{-\infty})$ means that
$\lim\limits_{\tau\rightarrow\infty}|O(\tau^{-\infty}) \tau^{s}|=0$ for any $s>0$.

Since $\beta>\beta_{p^{\ast}}$ and $A_{\tau}\geq1$, it follows from \eqref{eq2.1}--\eqref{a6} that
\begin{equation*}
d_{\beta}(p^{\ast})\leq E_{p^{\ast}}^{\beta}(u_{\tau})
=\frac{a}{2}A_{\tau}^{2}\tau^{2}+\frac{b}{4}\tau^{4}\left(A_{\tau}^{4}-\frac{\beta}{\beta_{p^{\ast}}}A_{\tau}^{p^{\ast}+2}\right)
+\frac{1}{2}V(x_{0})+o(1)+O(\tau^{-\infty})\rightarrow-\infty
\end{equation*}
as $\tau\rightarrow\infty$. This shows that $d_{\beta}(p^{\ast})=-\infty$ when  $\beta>\beta_{p^{\ast}}$.  So, there is no minimizer for $d_{\beta}(p^{\ast})$ if $\beta>\beta_{p^{\ast}}.$ $\hfill \square$\\

{\bf Proof of Theorem \ref{th3} :}
Let $u_{p} \geq 0$ be a  minimizer of \eqref{eq1.3}, it follows from \eqref{eq1.10} and \eqref{eq1.8} that
\begin{align} \label{hs}
\frac{a}{2}\int_{\mathbb R^{N}}|\nabla u_{p}|^{2}dx + \frac{b}{4}\left(\int_{\mathbb R^{N}}|\nabla u_{p}|^{2}dx\right)^{2} & + \frac{1}{2}\int_{\mathbb R^{N}}V(x)u_{p}^{2}dx
=d_{\beta}(p)+\frac{\beta}{p+2}\int_{\mathbb R^{N}}|u_{p}|^{p+2}dx \nonumber\\
&\leq  d_{\beta}(p) +\frac{b\beta}{4\beta_{p}}\left(\int_{\mathbb R^{N}}|\nabla u_{p}|^{2}dx\right)^{\frac{Np}{4}},
\end{align}
However, since $\beta>0$, by the definitions of $d_{\beta}(p)$ and $E_{p}^{\beta}$ it is clear that, for any $\xi(x)\in C_{0}^{\infty}(\mathbb R^{N})\cap S_1$, we have
\[
d_{\beta}(p)\leq E_{p}^{\beta}(\xi)\leq
\frac{a}{2}\int_{\mathbb R^{N}}|\nabla \xi|^{2}dx + \frac{b}{4}\left(\int_{\mathbb R^{N}}|\nabla \xi|^{2}dx\right)^{2}  + \frac{1}{2}\int_{\mathbb R^{N}}V(x)\xi^{2}dx:=
C \ (\text{independent of }p).
\]
Then,
\begin{equation}\label{c2}
\frac{b}{4}\left(\int_{\mathbb R^{N}}|\nabla u_{p}|^{2}dx\right)^{2}\leq C+\frac{b\beta}{4\beta_{p}}\left(\int_{\mathbb R^{N}}|\nabla u_{p}|^{2}dx\right)^{\frac{Np}{4}}.
\end{equation}
This implies that, there exists $M>0$ such that
\begin{equation}\label{a8}
\limsup_{p\nearrow p^{\ast}} \int_{\mathbb R^{N}}|\nabla u_{p}|^{2}dx \leq M.
\end{equation}
Otherwise, if $ \int_{\mathbb R^{N}}|\nabla u_{p}|^{2}dx  \rightarrow\infty $ as $p\nearrow p^{\ast}$,
by noting that $2=\frac{Np^*}{4}>\frac{Np}{4}$ and $\beta\in(0,\beta_{p^{\ast}})$, then, for  $p$ close to $p^{\ast}$,  it follows from  \eqref{c2}  that,
  \[
   \frac{b}{4} \leq \lim_{p\nearrow p^{\ast}} \frac{b}{4}\left(\int_{\mathbb R^{N}}|\nabla u_{p}|^{2}dx\right)^{\frac{N}{4}(p^*-p)} \leq   \frac{b\beta}{4\beta_{p^*}} <  \frac{b}{4},
  \]
  this leads to a contradiction.
Hence, \eqref{a8} holds and using again \eqref{hs} we know that
$\{u_{p}\}$ is bounded in $\mathcal{H}$. So, for some $u_{0}\in \mathcal{H}$, we may assume that
\begin{equation*}
u_{p}\rightharpoonup u_{0}\ \mbox{weakly in}\ \mathcal{H},\ \mbox{and}\ \ u_{p}\rightarrow u_{0}\ \mbox{strongly in}\ L^{q}(\mathbb R^{N}),\text{ with } q \in [2,2^{\ast}),
\end{equation*}
as $p\nearrow p^{\ast}$. By H\"{o}lder inequality,
\begin{equation*}
\int_{\mathbb R^{N}}|u_{p}|^{p+2}dx\leq\left(\int_{\mathbb R^{N}}|u_{p}|^{2}dx\right)^{\frac{p^{\ast}-p}{p^{\ast}}}
\left(\int_{\mathbb R^{N}}|u_{p}|^{p^{\ast}+2}dx\right)^{\frac{p}{p^{\ast}}}.
\end{equation*}
Note that $u_{p}$ is a minimizer of \eqref{eq1.3} and $p^*+2<2^*$ since $N<4$, then
\begin{equation}\label{a9}
\aligned
&\liminf_{p\nearrow p^{\ast}}d_{\beta}(p)=\liminf_{p\nearrow p^{\ast}}\Bigg\{\frac{a}{2}\int_{\mathbb R^{N}}|\nabla u_{p}|^{2}dx
+\frac{b}{4}\left(\int_{\mathbb R^{N}}|\nabla u_{p}|^{2}dx\right)^{2}
+\frac{1}{2}\int_{\mathbb R^{N}}V(x)u_{p}^{2}dx\\
&-\frac{\beta}{p+2}\int_{\mathbb R^{N}}|u_{p}|^{p+2}dx\Bigg\}\\
&\geq\liminf_{p\nearrow p^{\ast}}\Bigg\{\frac{a}{2}\int_{\mathbb R^{N}}|\nabla u_{p}|^{2}dx
+\frac{b}{4}\left(\int_{\mathbb R^{N}}|\nabla u_{p}|^{2}dx\right)^{2}
+\frac{1}{2}\int_{\mathbb R^{N}}V(x)u_{p}^{2}dx\\
&-\frac{\beta}{p+2}\left(\int_{\mathbb R^{N}}|u_{p}|^{p^{\ast}+2}dx\right)^{\frac{p}{p^{\ast}}}\Bigg\}\\
&\geq\frac{a}{2}\int_{\mathbb R^{N}}|\nabla u_{0}|^{2}dx
+\frac{b}{4}\left(\int_{\mathbb R^{N}}|\nabla u_{0}|^{2}dx\right)^{2}
+\frac{1}{2}\int_{\mathbb R^{N}}V(x)u_{0}^{2}dx
-\frac{\beta}{p^{\ast}+2}\int_{\mathbb R^{N}}|u_{0}|^{p^{\ast}+2}dx\\
&=E_{p^{\ast}}^{\beta}(u_{0})\geq d_{\beta}(p^{\ast}).
\endaligned
\end{equation}
But, by the definition of $ d_{\beta}$ we know that,  for any $\epsilon>0$, there exists  $u_{\epsilon}\in S_{1}$ such that
\begin{equation*}
E_{p^{\ast}}^{\beta}(u_{\epsilon})\leq d_{\beta}(p^{\ast})+\epsilon.
\end{equation*}
Therefore,
\begin{equation*}
\aligned
&\limsup_{p\nearrow p^{\ast}}d_{\beta}(p)\leq\limsup_{p\nearrow p^{\ast}}E_{p}^{\beta}(u_{\epsilon})\\
&=\limsup_{p\nearrow p^{\ast}}
\Bigg\{\frac{a}{2}\int_{\mathbb R^{N}}|\nabla u_{\epsilon}|^{2}dx
+\frac{b}{4}\left(\int_{\mathbb R^{N}}|\nabla u_{\epsilon}|^{2}dx\right)^{2}
+\frac{1}{2}\int_{\mathbb R^{N}}V(x)u_{\epsilon}^{2}dx
-\frac{\beta}{p+2}\int_{\mathbb R^{N}}|u_{\epsilon}|^{p+2}dx\Bigg\}\\
&=E_{p^{\ast}}^{\beta}(u_{\epsilon})
+\limsup_{p\nearrow p^{\ast}}\Bigg\{\frac{\beta}{p^{\ast}+2}\int_{\mathbb R^{N}}|u_{\epsilon}|^{p^{\ast}+2}dx
-\frac{\beta}{p+2}\int_{\mathbb R^{N}}|u_{\epsilon}|^{p+2}dx\Bigg\}\\
&\leq d_{\beta}(p^{\ast})+\epsilon.
\endaligned
\end{equation*}
So, combining \eqref{a9} and letting $\epsilon\rightarrow0$, we have
\begin{equation*}
\lim_{p\nearrow p^{\ast}}d_{\beta}(p)=d_{\beta}(p^{\ast})=E_{p^{\ast}}^{\beta}(u_{0}).
\end{equation*}
This shows that $u_{0}$ is a minimizer of $d_{\beta}(p^{\ast})$. Moreover,  $u_{p}\rightarrow u_{0}$ strongly in $\mathcal{H}$ as $p\nearrow p^{\ast}$. We finish the proof of Theorem \ref{th3}.
$\hfill \square$

Before going to prove the Theorem \ref{th2}, we need to establish some energy estimates on the minimizers of $\widetilde{d}_{\beta}(p)$.
By Lemma \ref{le6} and \ref{Le3.1}, we know that for any given $\beta>\beta_{p^{\ast}}$, $\widetilde{d}_{\beta}(p)$ has a minimizer when $p$ approaches $p^{\ast}$. Then, we have the following lemma.
\begin{lemma}\label{Le3.2}
For $1\leq N\leq3$, let $\beta>\beta_{p^{\ast}}$ be given and let $\widetilde{u}_{p} \geq 0$ be a   minimizer of $\widetilde{d}_{\beta}(p)$. Then
\begin{equation}\label{b7}
\left(\int_{\mathbb R^{N}}|\nabla \widetilde{u}_{p}|^{2}dx\right)^{2}
\approx\frac{4\beta}{b(p+2)}\int_{\mathbb R^{N}}|\widetilde{u}_{p}|^{p+2}dx
\approx\left(\frac{\beta p}{\beta_{p}p^{\ast}}\right)^{\frac{p^{\ast}}{p^{\ast}-p}}\triangleq r_{p},
\end{equation}
where, and in what follows,  we always use $m\approx n$ to denote that $\frac{m}{n}\rightarrow1$ as $p\nearrow p^{\ast}$.
\end{lemma}
Proof. The same as \eqref{c3}, we let $h(r)\triangleq\frac{b}{4}r-\frac{b\beta}{4\beta_{p}}r^{\frac{p}{p^{\ast}}}$.\\
We prove the lemma by contradiction. If there is a subsequence of $p\nearrow p^{\ast}$ such that, for some $\theta\geq0$,
\begin{equation*}
\frac{\left(\int_{\mathbb R^{N}}|\nabla \widetilde{u}_{p}|^{2}dx\right)^{2}}{r_{p}}\rightarrow \theta,  \mbox{ as }\ p\nearrow p^{\ast}.
\end{equation*}
We claim that there always a contradiction either $\theta \in [0,1)$, or $\theta >1$.

 In fact, if $\theta\in[0,1)$, then there exists $\epsilon>0$ such that $\delta\triangleq\theta+\epsilon<1$ and
$\frac{\left(\int_{\mathbb R^{N}}|\nabla \widetilde{u}_{p}|^{2}dx\right)^{2}}{r_{p}}\leq\delta$ as $p\nearrow p^{\ast}$. Hence,
\begin{equation*}
\lim_{p\nearrow p^{\ast}}\frac{h(\delta r_{p})}{h(r_{p})}
=\lim_{p\nearrow p^{\ast}}\frac{\delta r_{p}-\frac{\beta}{\beta_{p}}(\delta r_{p})^{\frac{p}{p^{\ast}}}}{r_{p}-\frac{\beta}{\beta_{p}}(r_{p})^{\frac{p}{p^{\ast}}}}
=\lim_{p\nearrow p^{\ast}}\frac{p^{\ast}\delta^{\frac{p}{p^{\ast}}}-p\delta}{p^{\ast}-p}
=\delta(1-\ln\delta)\in(0,1), \text{ for all }\ \delta\in(0,1).
\end{equation*}
For $p$ close to $p^{\ast}$, since $ {\left(\int_{\mathbb R^{N}}|\nabla \widetilde{u}_{p}|^{2}dx\right)^{2}}\leq \delta {r_{p}}$ and
$h(r)$ has a unique minimum point at $r_{p}$, these properties of the function $h$ imply that $0> \widetilde{d}_{\beta}(p)\geq h(\delta r_{p})$. Then
\begin{equation*}
\lim_{p\nearrow p^{\ast}}\frac{\widetilde{d}_{\beta}(p)}{h(r_{p})}
\leq\lim_{p\nearrow p^{\ast}}\frac{h(\delta r_{p})}{h(r_{p})}
=\delta(1-\ln\delta)<1,
\end{equation*}
this  contradicts  Lemma \ref{Le3.1}. \\
Similarly, if $\theta>1$, we have also a contradiction.

Thus, $\frac{\left(\int_{\mathbb R^{N}}|\nabla \widetilde{u}_{p}|^{2}dx\right)^{2}}{r_{p}}\rightarrow1$ as $p\nearrow p^{\ast}$.
Moreover, by \eqref{eq1.13}, $r_{p}\rightarrow+\infty$ as $p\nearrow p^{\ast}$, we then have
\begin{equation*}
\frac{\int_{\mathbb R^{N}}|\nabla \widetilde{u}_{p}|^{2}dx}{r_{p}}
=\frac{\left(\int_{\mathbb R^{N}}|\nabla \widetilde{u}_{p}|^{2}dx\right)^{2}}{r_{p}}\cdot\frac{1}{\int_{\mathbb R^{N}}|\nabla \widetilde{u}_{p}|^{2}dx}
\rightarrow0\quad\mbox{as}\ p\nearrow p^{\ast}.
\end{equation*}

Finally, since
\begin{equation*}
\frac{\widetilde{d}_{\beta}(p)}{r_{p}}
=\frac{a}{2r_{p}}\int_{\mathbb R^{N}}|\nabla \widetilde{u}_{p}|^{2}dx
+\frac{b}{4r_{p}}\left(\int_{\mathbb R^{N}}|\nabla \widetilde{u}_{p}|^{2}dx\right)^{2}
-\frac{\beta}{(p+2)r_{p}}\int_{\mathbb R^{N}}|\widetilde{u}_{p}|^{p+2}dx,
\end{equation*}
it follows from Lemma \ref{Le3.1} and  $p\nearrow p^{\ast}$ that
\begin{equation*}
\lim_{p\nearrow p^{\ast}}\frac{\beta}{(p+2)r_{p}}\int_{\mathbb R^{N}}|\widetilde{u}_{p}|^{p+2}dx=\frac{b}{4},
\end{equation*}
that is, $\frac{4\beta}{b(p+2)r_{p}}\int_{\mathbb R^{N}}|\widetilde{u}_{p}|^{p+2}dx\rightarrow1$ as $p\nearrow p^{\ast}$ and the lemma is proved.
$\quad\quad\quad\quad\quad\quad\quad\quad\square$

\begin{lemma}\label{le5}
For $1 \leq N\leq3$, let $\widetilde{u}_{p}$ be a nonnegative minimizer of $\widetilde{d}_{\beta}(p)$, set
\begin{equation}\label{eq5.6}
\epsilon_{p}=r_{p}^{-\frac{1}{4}} \text{ with } r_{p}=\left(\frac{\beta p}{\beta_{p}p^{\ast}}\right)^{\frac{p^{\ast}}{p^{\ast}-p}} \text{ as in \eqref{b7}},
 \mbox{ and } \widetilde{w}_{p}(x)=\epsilon_{p}^{\frac{N}{2}}\widetilde{u}_{p}(\epsilon_{p}x),
\end{equation}
then, there exist $w_{0}\in H^{1}(\mathbb R^{N})$ and a sequence $\{\widetilde{y}_{\epsilon_{p}}\}\subset\mathbb R^{N}$ such that, up to a subsequence,
\begin{equation*}
w_{p}(x)=\widetilde{w}_{p}(x+\widetilde{y}_{\epsilon_{p}})\rightarrow w_{0}\ \mbox{strongly in}\ H^{1}(\mathbb R^{N})\ \mbox{as}\ p\nearrow p^{\ast},
\end{equation*}
and
\begin{equation*}
w_{0}=\frac{1}{\|\phi_{p^{\ast}}\|_{L^{2}}}\phi_{p^{\ast}}(|x-x_{0}|)\quad \mbox{for}\ \mbox{some}\ x_{0}\in\mathbb R^{N}.
\end{equation*}
\end{lemma}
Proof. By Lemma \ref{Le3.1}, we have
\begin{equation}\label{c4}
\int_{\mathbb R^{N}}\widetilde{w}_{p}^{2}dx=\int_{\mathbb R^{N}}\widetilde{u}_{p}^{2}dx=1,
\end{equation}
\begin{equation}\label{c5}
\left(\int_{\mathbb R^{N}}|\nabla\widetilde{w}_{p}|^{2}dx\right)^{2}
=\epsilon_{p}^{4}\left(\int_{\mathbb R^{N}}|\nabla \widetilde{u}_{p}|^{2}dx\right)^{2}
\approx r_{p}^{-1}\cdot r_{p}=1\ \mbox{as}\ p\rightarrow p^{\ast}.
\end{equation}
\begin{equation}\label{eq3.9}
\int_{\mathbb R^{N}}|\widetilde{w}_{p}|^{p+2}dx
=\epsilon_{p}^{\frac{Np}{2}}\int_{\mathbb R^{N}}|\widetilde{u}_{p}|^{p+2}dx
\approx r_{p}^{-\frac{p}{p^{\ast}}}\cdot\frac{b(p+2)}{4\beta}r_{p}
\approx\frac{b(p^{\ast}+2)}{4\beta_{p^{\ast}}}\ \mbox{as}\ p\rightarrow p^{\ast}.
\end{equation}
By $\widetilde{u}_{p}$ is a minimizer of $\widetilde{d}_{\beta}(p)$, we know that  there exists $\lambda_{p}\in\mathbb R$ (Lagrange multiplier) such that
\begin{equation}\label{eq3.10}
-\left(a+b\int_{\mathbb R^{N}}|\nabla \widetilde{u}_{p}|^{2}dx\right)\triangle \widetilde{u}_{p}=\beta|\widetilde{u}_{p}|^{p}\widetilde{u}_{p}+\lambda_{p}\widetilde{u}_{p}.
\end{equation}
Then,
\begin{equation*}
\aligned
\lambda_{p}&=a\int_{\mathbb R^{N}}|\nabla \widetilde{u}_{p}|^{2}dx
+b\left(\int_{\mathbb R^{N}}|\nabla \widetilde{u}_{p}|^{2}dx\right)^{2}
-\beta\int_{\mathbb R^{N}}|\widetilde{u}_{p}|^{p+2}dx\\
&=4\widetilde{d}_{\beta}(p)
-a\int_{\mathbb R^{N}}|\nabla \widetilde{u}_{p}|^{2}dx
-\frac{\beta(p-2)}{p+2}\int_{\mathbb R^{N}}|\widetilde{u}_{p}|^{p+2}dx.
\endaligned
\end{equation*}
By Lemmas \ref{Le3.1}, \ref{Le3.2} and \eqref{eq5.6}, we see that
\begin{equation}\label{eq3.11}
\aligned
\lambda_{p}\epsilon_{p}^{4}
&=4\widetilde{d}_{\beta}(p)\cdot r_{p}^{-1}
-a\int_{\mathbb R^{N}}|\nabla \widetilde{u}_{p}|^{2}dx\cdot r_{p}^{-1}\\
&-\frac{b(p-2)}{4}\cdot\frac{4\beta}{b(p+2)}\int_{\mathbb R^{N}}|\widetilde{u}_{p}|^{p+2}dx\cdot r_{p}^{-1}
\rightarrow-\frac{b(4-N)}{2N}\quad\mbox{as}\ p\nearrow p^{\ast}.
\endaligned
\end{equation}
Since \eqref{eq3.9}, by  the concentration-compactness principle we know that there is a sequence $\{\widetilde{y}_{\epsilon_{p}}\}\subset\mathbb R^{N}$, and R$,\gamma>0$ such that
\begin{equation*}
\liminf_{p\nearrow p^{\ast}}\int_{B_{R}(\widetilde{y}_{\epsilon_{p}})}|\widetilde{w}_{p}|^{2}dx\geq\gamma>0.
\end{equation*}
Let
\begin{equation}\label{eq3.122}
w_{p}(x)=\epsilon_{p}^{\frac{N}{2}}\widetilde{u}_{p}(\epsilon_{p}x+\epsilon_{p}\widetilde{y}_{\epsilon_{p}})
=\widetilde{w}_{p}(x+\widetilde{y}_{\epsilon_{p}}),
\end{equation}
then
\begin{equation}\label{eq3.12}
\liminf_{p\nearrow p^{\ast}}\int_{B_{R}(0)}|w_{p}|^{2}dx\geq\gamma>0.
\end{equation}
By $\widetilde{u}_{p}$ satisfying \eqref{eq3.10}, we then know that $w_{p}(x)$ satisfies
\begin{equation}\label{eq3.13}
-a\epsilon_{p}^{2}\triangle w_{p}
-b\int_{\mathbb R^{N}}|\nabla w_{p}|^{2}dx\cdot\triangle w_{p}
=\lambda_{p}\epsilon_{p}^{4}w_{p}
+\beta\epsilon_{p}^{(4-\frac{Np}{2})}w_{p}^{p+1}.
\end{equation}
Using \eqref{eq3.7} and \eqref{eq5.6}, we see that
\begin{equation}\label{eq3.14}
\beta\epsilon_{p}^{(4-\frac{Np}{2})}
=\beta\epsilon_{p}^{\frac{N}{2}(p^{\ast}-p)}
=\beta\left(r_{p}^{\frac{p^{\ast}-p}{p^{\ast}}}\right)^{-1}
=\frac{\beta_{p}p^{\ast}}{p}\rightarrow \beta_{p^{\ast}}\quad\mbox{as}\ p\nearrow p^{\ast}.
\end{equation}
For any $\eta\in C_{0}^{\infty}(\mathbb R^{N})$, it follows from \eqref{eq1.13}, \eqref{eq3.7} and H\"{o}lder inequality that
\begin{equation*}
\left|\epsilon_{p}^{2}\int_{\mathbb R^{N}}\nabla w_{p}\nabla \eta dx\right|
\leq C\epsilon_{p}^{2}\left(\int_{\mathbb R^{N}}|\nabla w_{p}|^{2}dx\right)^{\frac{1}{2}}
=Cr_{p}^{-\frac{3}{4}}r_{p}^{\frac{1}{4}}(1+o(1))\rightarrow0\quad\mbox{as}\ p\nearrow p^{\ast}.
\end{equation*}
By \eqref{c4}, \eqref{c5} and \eqref{eq3.122}, we have
\begin{equation}\label{eq5.8}
\int_{\mathbb R^{N}}w_{p}^{2}dx=1\quad\mbox{and}\quad
\int_{\mathbb R^{N}}|\nabla w_{p}|^{2}dx
=\int_{\mathbb R^{N}}|\nabla \widetilde{w}_{p}|^{2}dx\rightarrow1\quad\mbox{as}\ p\nearrow p^{\ast}.
\end{equation}
Hence, $\{w_{p}\}$ is bounded in $H^{1}(\mathbb R^{N})$, we may assume that
\begin{equation}\label{c6}
w_{p}\rightharpoonup w_{0}\geq0\quad \mbox{weakly in}\ H^{1}(\mathbb R^{N})\quad\mbox{as}\ p\nearrow p^{\ast}
\end{equation}
for some $w_{0}\in H^{1}(\mathbb R^{N})$ and $w_{0}\not\equiv0$ by \eqref{eq3.12}. Moreover, it follows from \eqref{eq3.11}, \eqref{eq3.13}-\eqref{c6} that $w_{0}$ satisfies
\begin{equation}\label{eq3.15}
-b\triangle w_{0}=-\frac{b(4-N)}{2N}w_{0}+\beta_{p^{\ast}}w_{0}^{p^{\ast}+1}.
\end{equation}
Combining \eqref{eq3.15} and the Pohozaev identity, we have
\begin{equation}\label{eq5.7}
\begin{cases}
\int_{\mathbb R^{N}}|\nabla w_{0}|^{2}dx=\int_{\mathbb R^{N}}w_{0}^{2}dx,
\\\\ \int_{\mathbb R^{N}}|w_{0}|^{p^{\ast}+2}dx=\frac{b(4+N)}{2N\beta_{p^{\ast}}}\int_{\mathbb R^{N}}w_{0}^{2}dx.
\end{cases}
\end{equation}
Then, it follows from \eqref{eq1.8} and \eqref{eq5.7} that
\begin{equation}\label{eq5.9}
\frac{2N\beta_{p^{\ast}}}{b(4+N)}\leq\frac{\left(\int_{\mathbb R^{N}}|\nabla w_{0}|^{2}dx\right)^{2}\left(\int_{\mathbb R^{N}}w_{0}^{2}dx\right)^{\frac{4-N}{N}}}{\int_{\mathbb R^{N}}|w_{0}|^{p^{\ast}+2}dx}
=\frac{2N\beta_{p^{\ast}}}{b(4+N)}\left(\int_{\mathbb R^{N}}w_{0}^{2}dx\right)^{\frac{4}{N}},
\end{equation}
this shows that
\begin{equation*}
\int_{\mathbb R^{N}}w_{0}^{2}dx\geq1,
\end{equation*}
which together with \eqref{eq5.8} and \eqref{c6}, we have
\begin{equation}\label{eq5.10}
\int_{\mathbb R^{N}}w_{0}^{2}dx=1.
\end{equation}
Hence,
\begin{equation*}
w_{p}\rightarrow w_{0}\ \mbox{strongly in}\ L^{2}(\mathbb R^{N})\quad\mbox{as}\ p\nearrow p^{\ast}.
\end{equation*}
and it follows from  \eqref{eq5.8} and \eqref{eq5.7} that
\begin{equation*}
\int_{\mathbb R^{N}}|\nabla w_{p}|^{2}dx\rightarrow\int_{\mathbb R^{N}}|\nabla w_{0}|^{2}dx\quad\mbox{as}\ p\nearrow p^{\ast}.
\end{equation*}
So, $w_{p}\rightarrow w_{0}\ \mbox{strongly in}\ H^{1}(\mathbb R^{N})$ as $p\nearrow p^{\ast}$.

Furthermore, \eqref{c6} and the strong maximum principle imply that $w_{0}>0$. Note that $w_{0} $ is a positive solution of \eqref{c7} and also of \eqref{eq1.9}  (up to a rescaling),  then the uniqueness of positive  solution of \eqref{eq1.9} implies that
\begin{equation}\label{eq5.11}
w_{0}=\frac{1}{\|\phi_{p^{\ast}}\|_{L^{2}}}\phi_{p^{\ast}}(|x-x_{0}|),\quad \mbox{for}\ \mbox{some}\ x_{0}\in\mathbb R^{N}.
\end{equation}
$\hfill \square$

Now, we turn to showing the decay property for $w_{p}$ defined by \eqref{eq3.122}. By Lemma \ref{le5}, we see that if $\widetilde{u}_{p}$ is a nonnegative minimizer of $\widetilde{d}_{\beta}(p)$, then
there exist a subsequence $\{p_{k}\}$ with $p_{k}\nearrow p^{\ast}$ as $k\rightarrow\infty$ and a positive function $w_{0}$ such that
\begin{equation}\label{eq3.16}
w_{p_{k}}\rightarrow w_{0}\quad \mbox{strongly in}\ H^{1}(\mathbb R^{N})\quad \mbox{as}\ k\rightarrow\infty,
\end{equation}
Hence, for any $\alpha\in[2,2^{\ast})$,
\begin{equation}\label{eq5.1}
\int_{|x|\geq R}|w_{p_{k}}|^{\alpha}dx\rightarrow0\quad\mbox{as}\  R\rightarrow\infty\  \mbox{uniformly}\  \mbox{for large}\ k.
\end{equation}
By \eqref{eq3.11}, \eqref{eq3.13} and \eqref{eq3.14},  we know that
\begin{equation*}
-\triangle w_{p_{k}}-c(x)w_{p_{k}}\leq0,\quad\mbox{for large}\ k,
\end{equation*}
where $c(x)=\frac{2\beta_{p^{\ast}}}{b}w_{p}^{p}.$ By applying De-Giorgi-Nash-Morse theory (similar to the proof of \cite[Theorem 4.1]{HL2}), we deduce that
\begin{equation}\label{eq5.2}
\max_{B_{1}(\xi)}w_{p_{k}}\leq C\left(\int_{B_{2}(\xi)}|w_{p_{k}}|^{2}dx\right)^{\frac{1}{2}},
\end{equation}
where $\xi$ is an arbitrary point in $\mathbb R^{N}$, and $C$ is a constant depending only on the bound of $\|w_{p_{k}}\|_{L^{3}(B_{2}(\xi))}$.
Hence,
\begin{equation}\label{eq5.3}
w_{p_{k}}(x)\rightarrow0 \ \mbox{as}\ |x|\rightarrow\infty, \quad \mbox{uniformly for large}\ k.
\end{equation}
Note from \eqref{eq3.11},\eqref{eq3.13}, \eqref{eq3.14} and \eqref{eq5.3} that there exists a $R>0$, independent of $k$, such that $w_{p_{k}}$ satisfies
\begin{equation}\label{eq5.4}
-\Delta w_{p_{k}}(x)+\frac{4-N}{4N}w_{p_{k}}(x)\leq0\quad\mbox{uniformly for large}\ k\ \mbox{and}\ |x|>R.
\end{equation}
Applying the comparison principle \cite{KW} to compare $w_{p_{k}}$ with $Ce^{-\sqrt{\frac{4-N}{4N}}|x|}$, we then know that there exists $C>0$, independent of $k$, such that
\begin{equation}\label{eq5.5}
w_{p_{k}}(x)\leq Ce^{-\sqrt{\frac{4-N}{4N}}|x|}\quad\mbox{uniformly for large}\ k\ \mbox{and}\ |x|>R.
\end{equation}
\begin{lemma}\label{Le3.3}
If $1 \leq N\leq3$ and
 $V(x)$ satisfies \eqref{a3}, let $\beta>\beta_{p^{\ast}}$  and $u_{p}$ be a  minimizer of \eqref{eq1.3}. Then
\begin{equation}\label{eq3.17}
0\leq d_{\beta}(p)-\widetilde{d}_{\beta}(p)\rightarrow0\quad as\ p\nearrow p^{\ast},
\end{equation}
and
\begin{equation}\label{eq3.18}
\int_{\mathbb R^{N}}V(x)u_{p}^{2}dx\rightarrow0\quad as\ p\nearrow p^{\ast}.
\end{equation}
\end{lemma}
Proof. By the definition of $d_{\beta}(p)$ and $\widetilde{d}_{\beta}(p),$ it is easy to see that
\begin{equation*}
d_{\beta}(p)-\widetilde{d}_{\beta}(p)\geq0.
\end{equation*}
Now we come to get an upper bound for $d_{\beta}(p)-\widetilde{d}_{\beta}(p)$. Let $\xi(x)\in C_{0}^{\infty}(\mathbb R^{N})$ and $0\leq\xi(x)\leq1$ such that $\xi(x)\equiv1$ if $|x|\leq1$, $\xi(x)\equiv0$ if $|x|\geq2$, and $|\nabla\xi(x)|\leq C_{0}$. For any $x_{0}\in\mathbb{R}^{N}$, we take
\begin{equation*}
\overline{u}_{p}(x)
=A_{p}\xi(x-x_{0})\epsilon_{p}^{-\frac{N}{2}}w_{p}\left(\frac{x-x_{0}}{\epsilon_{p}}\right)
=A_{p}\xi(x-x_{0})\widetilde{u}_{p}(x-x_{0}+\epsilon_{p}\widetilde{y}_{\epsilon_{p}}),
\end{equation*}
where $A_{p}>0$ is chosen so that $\|\overline{u}_{p}\|_{L^{2}}^{2}=1$ and $w_{p}$ is given by \eqref{eq3.122}. Then, it follows from \eqref{eq5.8} that
$\int_{\mathbb{R}^{N}}|w_{p}|^{2}dx=1$ and
\begin{equation*}
1\leq A_{p}^{2}=
\frac{\int_{\mathbb R^{N}}|w_{p}|^{2}dx}{\int_{\mathbb R^{N}}\xi^{2}(\epsilon_{p}x)w_{p}^{2}(x)dx}
\leq\frac{\int_{\mathbb R^{N}}|w_{p}|^{2}dx}{\int_{|\epsilon_{p}x|\leq1}\xi^{2}(\epsilon_{p}x)w_{p}^{2}(x)dx}.
\end{equation*}
By \eqref{eq1.13} and \eqref{eq5.5}, we see that
\begin{equation}\label{eqa}
0\leq A_{p}^{2}-1\leq\frac{\int_{|\epsilon_{p}x|\geq1}|w_{p}|^{2}dx}{\int_{|\epsilon_{p}x|\leq1}\xi^{2}(\epsilon_{p}x)w_{p}^{2}(x)dx}
\leq Ce^{-\sqrt{\frac{4-N}{4N}}\epsilon_{p}^{-1}}
\end{equation}
as $p\nearrow p^{\ast}$, and
\begin{equation}\label{eqb}
1\leq A_{p}^{p+2}\leq
\left(1+Ce^{-\sqrt{\frac{4-N}{4N}}\epsilon_{p}^{-1}}\right)^{\frac{p+2}{2}}
\leq1+6Ce^{-\sqrt{\frac{4-N}{4N}}\epsilon_{p}^{-1}}.
\end{equation}
By \eqref{eqa},
\begin{equation}\label{eqc}
\int_{\mathbb R^{N}}V(x)\overline{u}_{p}^{2}dx
=A_{p}^{2}\int_{\mathbb R^{N}}V(\epsilon_{p}x+x_{0})\xi^{2}(\epsilon_{p}x)w_{p}^{2}(x)dx
\rightarrow V(x_{0})\int_{\mathbb R^{N}}w_{0}^{2}dx=V(x_{0})
\end{equation}
as $p\nearrow p^{\ast}$, and
\begin{equation}\label{eqd}
\aligned
\int_{\mathbb R^{N}}|\overline{u}_{p}|^{p+2}
&=\epsilon_{p}^{-\frac{NP}{2}}A_{p}^{p+2}\int_{\mathbb R^{N}}\xi^{p+2}(\epsilon_{p}x)w_{p}^{p+2}(x)dx\\
&=(A_{p}^{p+2}-1)\epsilon_{p}^{-\frac{NP}{2}}\int_{\mathbb R^{N}}\xi^{p+2}(\epsilon_{p}x)w_{p}^{p+2}(x)dx\\
&+\epsilon_{p}^{-\frac{NP}{2}}\int_{\mathbb R^{N}}(\xi^{p+2}(\epsilon_{p}x)-1)w_{p}^{p+2}(x)dx
+\epsilon_{p}^{-\frac{NP}{2}}\int_{\mathbb R^{N}}w_{p}^{p+2}(x)dx\\
&\leq\int_{\mathbb R^{N}}\widetilde{u}_{p}^{p+2}(x)dx
+ Ce^{-\sqrt{\frac{4-N}{4N}}\epsilon_{p}^{-1}}
\endaligned
\end{equation}
as $p\nearrow p^{\ast}$. Similarly, we know that
\begin{equation}\label{eqe}
\int_{\mathbb R^{N}}|\nabla\overline{u}_{p}|^{2}dx
\leq\int_{\mathbb R^{N}}|\nabla\widetilde{u}_{p}|^{2}dx+ Ce^{-\sqrt{\frac{4-N}{4N}}\epsilon_{p}^{-1}}
\end{equation}
and
\begin{equation}\label{eqf}
\left(\int_{\mathbb R^{N}}|\nabla\overline{u}_{p}|^{2}dx\right)^{2}
\leq\left(\int_{\mathbb R^{N}}|\nabla\widetilde{u}_{p}|^{2}dx\right)^{2}+ Ce^{-\sqrt{\frac{4-N}{4N}}\epsilon_{p}^{-1}}
\end{equation}
as $p\nearrow p^{\ast}$. Taking $x_{0}\in\mathbb R^{N}$ such that $V(x_{0})=0$, then we deduce from \eqref{eqc}-\eqref{eqf} that
\begin{equation*}
\aligned
0\leq d_{\beta}(p)-\widetilde{d}_{\beta}(p)
&\leq E_{p}^{\beta}(\overline{u}_{p})-\widetilde{E}_{p}^{\beta}(\widetilde{u}_{p})\\
&=\widetilde{E}_{p}^{\beta}(\overline{u}_{p})-\widetilde{E}_{p}^{\beta}(\widetilde{u}_{p})+\frac{1}{2}\int_{\mathbb R^{N}}V(x)\overline{u}_{p}^{2}dx\\
&\leq\frac{1}{2}V(x_{0})+ Ce^{-\sqrt{\frac{4-N}{4N}}\epsilon_{p}^{-1}}+o(1)\rightarrow0
\endaligned
\end{equation*}
as $p\nearrow p^{\ast}$. Furthermore, if $u_{p}$ is a minimizer of $d_{\beta}(p)$,  we should have
\begin{equation*}
0\leq\frac{1}{2}\int_{\mathbb R^{N}}V(x)u_{p}^{2}dx=d_{\beta}(p)-\widetilde{E}_{p}^{\beta}(u_{p}(x))
\leq d_{\beta}(p)-\widetilde{d}_{\beta}(p)\rightarrow0\quad\mbox{as}\ p\nearrow p^{\ast}.
\end{equation*}
$\quad\quad\quad\quad\quad\quad\quad\quad\quad\quad\quad\quad\quad\quad\quad\quad\quad\quad\quad\quad\quad\quad\quad\quad\quad\quad\quad\quad
\quad\quad\quad\quad\quad\quad\quad\quad\quad\quad\quad\quad\quad\quad\quad\quad\quad\square$\\
{\bf Proof of Theorem \ref{th2} :} Let $u_{p} \geq 0$ be a   minimizer of $d_{\beta}(p)$.  Using \eqref{eq3.17} and \eqref{eq3.18}, similar to the proof of Lemma \ref{Le3.2}, we have also that
\begin{equation}\label{eq3.19}
\left(\int_{\mathbb R^{3}}|\nabla u_{p}|^{2}dx\right)^{2}
\approx\frac{4\beta}{b(p+2)}\int_{\mathbb R^{N}}|u_{p}|^{p+2}dx
\approx\left(\frac{\beta p}{\beta_{p}p^{\ast}}\right)^{\frac{p^{\ast}}{p^{\ast}-p}}=r_{p}
\end{equation}
Similar to \eqref{eq3.12}, there exist  $\{\overline{y}_{\epsilon_{p}}\}\subset\mathbb R^{N}$, and $R,\gamma>0$ such that
\begin{equation}\label{c7}
\liminf_{p\nearrow p^{\ast}}\int_{B_{R}(0)}|\overline{w}_{p}|^{2}dx \geq \gamma>0,
\end{equation}
where
\begin{equation*}
\overline{w}_{p}(x)=\epsilon_{p}^{\frac{N}{2}}u_{p}(\epsilon_{p}x+\epsilon_{p}\overline{y}_{\epsilon_{p}}),  \text{ where } \epsilon_{p} \text{ is  given by } \eqref{eq5.6}.
\end{equation*}
Since $u_{p}$ is a minimizer of $d_{\beta}(p)$, similar to \eqref{eq3.13}, there exists $\overline{\lambda}_{p}$ such that
\begin{equation}\label{x6}
-a\epsilon_{p}^{2}\triangle \overline{w}_{p}
-b\int_{\mathbb R^{N}}|\nabla \overline{w}_{p}|^{2}dx\cdot\triangle \overline{w}_{p}
+\epsilon_{p}^{4}V(\epsilon_{p}x+\epsilon_{p}\overline{y}_{\epsilon_{p}})\overline{w}_{p}
=\overline{\lambda}_{p}\epsilon_{p}^{4}\overline{w}_{p}
+\beta\epsilon_{p}^{(4-\frac{Np}{2})}\overline{w}_{p}^{p+1}.
\end{equation}
We claim that $\{\epsilon_{p}\overline{y}_{\epsilon_{p}}\}$ is bounded uniformly in $p\nearrow p^{\ast}$.
Otherwise,
we may assume that, there is a subsequence  $p_{n}\nearrow p^{\ast}$ as $n\rightarrow\infty$ such that
\begin{equation*}
\epsilon_{p_{n}}|\overline{y}_{\epsilon_{p_{n}}}|\rightarrow\infty
\text{ and }
\epsilon_{p_{n}}\rightarrow0 \mbox{ as }\ n\rightarrow\infty.
\end{equation*}
Since \eqref{eq3.18},  we have
\begin{equation}\label{eqg}
\int_{\mathbb R^{N}}V(x)u^{2}_{p}(x)dx
=\int_{\mathbb R^{N}}V(\epsilon_{p}x+\epsilon_{p}\overline{y}_{\epsilon_{p}})\overline{w}_{p}^{2}(x)dx\rightarrow0\quad\mbox{as}\ p\nearrow p^{\ast}.
\end{equation}
By \eqref{a3}, there exists $C_{0}>0$ such that $V(x)\geq C_{0}$ for $|x|$ being large enough. We then derive from \eqref{c7}
and Fatou's Lemma that
\begin{equation*}
\liminf_{n\rightarrow\infty}\int_{\mathbb R^{N}}V(\epsilon_{p_{n}}x+\epsilon_{p_{n}}\overline{y}_{\epsilon_{p_{n}}})\overline{w}_{p_{n}}^{2}(x)dx
\geq\int_{\mathbb R^{N}}\liminf_{n\rightarrow\infty}V(\epsilon_{p_{n}}x+\epsilon_{p_{n}}\overline{y}_{\epsilon_{p_{n}}})\overline{w}_{p_{n}}^{2}(x)dx
\geq\gamma C_{0}>0,
\end{equation*}
which contradicts \eqref{eqg}.

So,  $\{\epsilon_{p}\overline{y}_{\epsilon_{p}}\}$ is bounded uniformly in $p\nearrow p^{\ast}$. By passing to a subsequence, we may assume that $\epsilon_{p}\overline{y}_{\epsilon_{p}}\rightarrow z_{0}$ as $p\nearrow p^{\ast}$ for some $z_{0}\in \mathbb R^{N}$. It follows from \eqref{c7} and Fatou's Lemma that
\begin{equation*}
\liminf_{p\nearrow p^{\ast}}\int_{\mathbb R^{N}}V(\epsilon_{p}x+\epsilon_{p}\overline{y}_{\epsilon_{p}})\overline{w}_{p}^{2}(x)dx
\geq V(z_{0})\int_{B_{R}(0)}\liminf_{p\nearrow p^{\ast}}\overline{w}_{p}^{2}(x)dx
\geq V(z_{0})\eta,
\end{equation*}
this together with \eqref{eqg} imply that $V(z_{0})=0$.

Finally, by \eqref{eq3.19} and similar to Lemma \ref{le5}, we know that $\overline{w}_{p}=\epsilon_{p}^{\frac{N}{2}}u_{p}(\epsilon_{p}x+\epsilon_{p}\overline{y}_{\epsilon_{p}})$ satisfies
\begin{equation*}
\overline{w}_{p}\rightarrow \overline{w}_{0}\quad\mbox{in}\ H^{1}(\mathbb R^{N})\quad\mbox{as}\ p\nearrow p^{\ast}.
\end{equation*}
and $\overline{w}_{0}$ satisfies the following equation
\begin{equation*}
-b\triangle \overline{w}_{0}=-\frac{b(4-N)}{2N}\overline{w}_{0}+\beta_{p^{\ast}}\overline{w}_{0}^{p^{\ast}+1}.
\end{equation*}
Thus,   the uniqueness (up to translations) of positive  solution of \eqref{eq1.9} implies that
\begin{equation*}
\overline{w}_{0}=\frac{1}{\|\phi_{p^{\ast}}\|_{L^{2}}}\phi_{p^{\ast}}(|x-y_{0}|)\quad \mbox{for}\ \mbox{some}\ y_{0}\in\mathbb R^{N}.
\end{equation*}
The proof is completed.$\quad\quad\quad\quad\quad\quad\quad\quad\quad\quad\quad\quad\quad\quad\quad\quad\quad\quad\quad\quad\quad\quad
\quad\quad\quad\quad\quad\quad\quad\quad\quad\quad\quad\quad\quad\square$

\end{document}